\documentclass[11pt]{article}
\usepackage{amsmath,amsthm,amsfonts,amssymb,amscd,dsfont}

\newtheorem{theorem}{Theorem}[section]
\newtheorem{lemma}[theorem]{Lemma}
\newtheorem{proposition}[theorem]{Proposition}
\newtheorem{corollary}[theorem]{Corollary}
\theoremstyle{definition}
\newtheorem{definition}[theorem]{Definition}
\newtheorem{algorithm}[theorem]{Algorithm}
\newtheorem{example}[theorem]{Example}

\newtheorem{remark}{Remark}[section]
\numberwithin{equation}{section}

\def\re{\mathbb R}

\def\lV{\left\Vert }
\def\rV{\right\Vert }

\begin{document}

\title{\Large\bf{Convergence analysis of the extragradient method for vector quasi-equilibrium problems}}

\author{
Vahid Mohebbi%
\thanks{Department of Mathematical Sciences, University of Texas at El Paso, 500 W. University Avenue, El Paso, Texas 79968, USA
({\tt vmohebbi@utep.edu; mohebbimath@gmail.com}).}
}

\maketitle

\begin{abstract}
We study the extragradient method for solving vector quasi-equilibrium problems in Banach spaces, which generalizes the extragradient method for vector equilibrium problems and scalar quasi-equilibrium problems.
We propose a regularization procedure which ensures strong convergence of the generated sequence to a solution of the vector quasi-equilibrium problem, under standard assumptions on the problem without assuming neither any monotonicity assumption on the vector valued bifunction nor any weak continuity assumption of $f$ in its arguments that in the many well-known methods have been used. Also, we show that the boundedness of the generated sequences implies that the solution set of the vector quasi-equilibrium problem is nonempty, and prove the strong convergence of the generated sequences to a solution of the problem. Finally, we give some examples of vector quasi-equilibrium problems in several Banach spaces to which  our main theorem can be applied. We also present some numerical experiments. 
\end{abstract}


\medskip

{\bf AMS Classification Number:} 90C25, 90C30.

\section{Introduction}\label{s0}

Let $E$ be a reflexive Banach space and $K\subset E$ be a nonempty, closed and convex set.  Suppose that $T(\cdot)$ is a multivalued mapping with nonempty values from $K$ into ${\cal P}(K)$  and  $f:E\times E\to\re$ is a bifunction.
The quasi-equilibrium problem QEP$(f,T)$ consists of finding $x^*\in T(x^*)$, i.e. a fixed point $x^*$ of $T(\cdot)$, such that
\begin{equation}\label{QEP}
f(x^*,y)\ge 0, \ \ \ \ \ \forall y\in T(x^*).
\end{equation}
 
The associated  Minty quasi-equilibrium problem can be expressed as finding $x^*\in T(x^*)$ such that $f(y,x^*)\leq 0$ for all $y\in T(x^*)$. 
When the constraint set $T(x)$ is equal to $K$ for every $x\in K$, the quasi-equilibrium problem QEP$(f,T)$ becomes a classical  equilibrium problem EP$(f,K)$, and the associated  Minty quasi-equilibrium problem becomes a classical  Minty equilibrium problem.\\
A simple example of  quasi-equilibrium problems is a quasi-variational inequality
problem. Let  $T(\cdot)$ be a multivalued mapping with nonempty values from $K$ into itself  and consider a map $A:E\to E^*$ where $E^*$ is the topological dual of $E$, and define $f(x,y)=\langle A(x),y-x\rangle$,
where $\langle\cdot ,\cdot\rangle:E^*\times E\to\re$ denotes the duality pair,
i.e. $\langle z,x\rangle=z(x)$. 
Then QEP$(f,T)$ is equivalent to the quasi-variational inequality
problem QVIP$(A,T)$, consisting of finding a point $x^*\in T(x^*)$ such that  $\langle
A(x^*),y-x^*\rangle\ge 0$ for all $y\in T(x^*)$. 

Equilibrium problems and quasi-equilibrium problem have 
been studied extensively in Hilbert, Hadamard, Banach as well as in 
topological vector spaces by many authors (e.g. \cite{ACI}, \cite{BiS},\cite{CCR}, \cite{CoH}, \cite{DJ-M1}, \cite{IKS},
 \cite{IuM}, \cite{ISo}, \cite{KM2}, \cite{KM1}, \cite{VSNV}).
 
Now we move from scalar valued bifunctions to vector valued ones. 
We assume that $K\subset E$ is a nonempty, closed and convex set, $T(\cdot)$ is a multivalued mapping with nonempty values from $K$ into ${\cal P}(K)$  and suppose that $Y$ is a
real Banach space containing a closed, convex and pointed cone $C$  with nonempty interior (denoted
as int$(C)$), and
$f:E\times E\rightarrow Y$ is a vector valued bifunction. The {\it vector equilibrium problem}
denoted as VEP$(f,K)$, consists of finding $x^*\in K$ such that
\begin{equation}
\label{VEP}
f(x^*,y)\not\in -{\rm int}(C),  \ \ \ \ \ \ \ \forall  y\in K.
\end{equation}
If $x^*$ satisfies \eqref{VEP}, then $x^*$ is said to be a solution or equilibrium point for VEP$(f,K)$. 
In addition, it is valuable to mention that many authors generalized the
equilibrium problem EP$(f,K)$ to the vector case in the following ways (see \cite{Gong}, \cite{KAS}):\\
find $x^*\in K$ such that
\begin{equation}\label{VEP-12}
f(x^*,y)\not\in -C\setminus \{0\}, \ \ \ \ \ \ \ \forall  y\in K,
\end{equation}
or find $x^*\in K$ such that
\begin{equation}\label{VEP-13}
f(x^*,y)\in C, \ \ \ \ \ \ \ \forall  y\in K.
\end{equation}
The associated dual vector equilibrium problem can be expressed as finding $x^*\in K$ such that
 \begin{equation}\label{DVEP}
 f(y,x^*)\in -C, \ \ \ \ \ \ \ \forall y\in K.
 \end{equation}

 The {\it vector quasi-equilibrium problem}
denoted as VQEP$(f,T)$, consists of finding $x^*\in T(x^*)$ such that
\begin{equation}
\label{QVEP}
f(x^*,y)\not\in -{\rm int}(C),  \ \ \ \ \ \ \ \forall  y\in T(x^*).
\end{equation}
If $x^*$ satisfies \eqref{QVEP}, then $x^*$ is said to be a solution or equilibrium point for VQEP$(f,T)$. We
denote the set of all equilibrium points of VQEP$(f,T)$ as $S(f,T)$.
We also denote the set of all fixed points of the multivalued mapping $T(\cdot)$ by ${\rm Fix}(T)$.
The associated dual vector quasi-equilibrium problem can be expressed as finding $x^*\in T(x^*)$ such that
 \begin{equation}\label{DQVEP}
 f(y,x^*)\in -C, \ \ \ \ \forall y\in T(x^*).
 \end{equation}

The prototypical example of vector equilibrium problems occurs when $Y=\re^m$ and $C$ is the nonnegative cone, i.e. 
$C=\re^m_+$. If we take $G:E\to\re^m$ and $f(x,y)=G(y)-G(x)$ then VEP$(f,K)$ is equivalent to the problem of finding
a Pareto minimizer of $G$ on $K$, i.e. a point $x^*\in K$ such that there exists no $x\in K$ such that
$G(x)\le G(x^*)$ and $G(x)\ne G(x^*)$ (here $G(x)\le G(x^*)$ means $G(x)_i\le G(x^*)_i$ for all $i\in \{1,\dots m\}$).

We will deal in this paper with the extragradient (or Korpelevich's) method 
for vector quasi-equilibrium problems in infinite dimensional Banach
spaces, and thus we
start with an introduction to its well known finite
dimensional formulation when applied to variational inequalities, i.e.,
we assume that $E=\re^n, Y=\re$ and $f(x,y)=\langle A(x),y-x\rangle$ with $A:\re^n\to\re^n$.
We assume that $A$ is {\it monotone}, i.e. $\langle A(x)-A(y),x-y\rangle\ge 0$ for all $x,y\in E$.
In this setting, there are several iterative methods
for solving VIP($A,K$). One of the most useful ones is the {\it extragradient method} presented in
\cite{Kor}, which generates a sequence $\{x^k\}\subset E$ according to: 
\begin{equation}\label{e2}
y^k=P_K(x^k-\alpha_kA(x^k)),
\end{equation}
\begin{equation}\label{e3}
x^{k+1}=P_K(x^k-\alpha_kA(y^k)),
\end{equation}
where $P_ K$ denotes the orthogonal projection onto $K$  and $\{\alpha_k\}\subset\re$ is a
sequence of positive stepsizes.

It was proved in \cite{Kor} that if $A$ is monotone and
Lipschitz continuous with constant $L$, and VIP($A,K$) has solutions,
then the sequence generated by \eqref{e2}--\eqref{e3} converges to a solution of
VIP($A,K$) provided that  $\alpha_k=\alpha \in (0, 1/L)$.

In the absence of Lipschitz continuity of $A$, it is natural to 
search for an
appropriate stepsize in an inner loop.
This is achieved in the following procedure:

Take
$\delta\in(0,1)$, $\hat{\beta}$, $\tilde{\beta}$ satisfying
$0<\hat{\beta}\leq\tilde{\beta}$, and a sequence
$\{\beta_k\}\subset[\hat{\beta},\tilde{\beta}]$. 
The method is initialized with any 
$x^0\in K$ and the iterative step is as follows:

Given $x^k$, define
\begin{equation}\label{paso1}
z^k:= x^k-\beta_k A(x^k) .
\end{equation}
\noindent If $x^k=P_K(z^k)$ stop. Otherwise take
\begin{equation}\label{j(k)}
j(k):=\min\left\{\,\,j\geq0:\left\langle
A(2^{-j}P_K(z^k)+(1-2^{-j})x^k),x^k-P_K(z^k)\right\rangle\geq
\frac{\delta}{\beta_k}\,\|x^k-P_K(z^k)\|^2\,\,\right\},
\end{equation}
\begin{equation}\label{alpha_k}
\alpha_k:=2^{-j(k)},
\end{equation}
\begin{equation}\label{armijo}
y^k:=\alpha_k P_K(z_k) + (1-\alpha_k)x^k,
\end{equation}
\begin{equation}\label{Hk2}
H_k:=\left\{z\in \re^n\,:\,\langle A(y^k),z-y^k\rangle =0\right\},
\end{equation}
\begin{equation}\label{paso3}
x^{k+1}:=P_K\left( P_{H_k}( x^k)\right).
\end{equation}
 
This method converges to a solution of VIP($A,K$) under the only assumptions of monotonicity
of $A$ and existence of solutions; see \cite{IS1}.

The above backtracking procedure for determining the right  $\alpha_k$ is sometimes
called an Armijo-type search (see \cite{Arm}). It has been analyzed for VIP($A,K$)
in \cite{IS1} and \cite{Ko1}. Other variants of Korpelevich's method can be found in \cite{GoT}, \cite{IuM}, \cite{KM3}, \cite{Kho},
\cite{Mar}, \cite{VSNV}. 

Recently an extragradient method for vector equilibrium problems in a Banach space $E$ 
has been studied in \cite{IM-3}. 
It has the following form: 

\begin{algorithm}\label{ka-19}
\end{algorithm}
Take
$\delta\in(0,1)$, $\hat\beta$, $\tilde\beta$ satisfying
$0<\hat{\beta}\leq\tilde{\beta}$, a sequence
$\{\beta_k\}\subset[\hat\beta,\tilde\beta]$, and a sequence $\{e^k\}\subset {\rm int}(C)$ such that $\lV e^k\rV=1$.
 
\noindent{\bf 1. Initialization}:

\begin{equation}\label{ka-1}
x^0 \in K.
\end{equation}

\noindent{\bf 2. Iterative step}:  Given $x^k$, define
\begin{equation}\label{ka-2}
 z^k \in {\rm argmin^C_w}\left\{f(x^k,y)+
\frac{1}{2\beta_k}\lV y\rV^2e^k-\frac{1}{\beta_k}\langle y, Jx^k \rangle e^k : y\in K \right\}.
\end{equation}
If $x^k=z^k$ stop. Otherwise, let
\begin{equation}\label{ka-4}
\ell(k)=\min\left\{\ell \geq 0: -\beta_kf(y^{\ell},x^k)+\beta_kf(y^{\ell},z^k)+ 
\frac{\delta}{2}\phi(z^k,x^k)e^k\not\in {\rm int}(C) \right\},
\end{equation}
where
\begin{equation}\label{ka-5}
y^{\ell}=2^{-\ell}z^k+(1-2^{-\ell})x^k.
\end{equation}
We take
\begin{equation}\label{ka-6}
\alpha_k:=2^{-\ell(k)},
\end{equation}
\begin{equation}\label{ka-7}
y^k:=\alpha_k z^k+(1-\alpha_k)x^k,
\end{equation}
\begin{equation}\label{ka-8}
w^k=P_{H_k}(x^k),
\end{equation}
where
$$
H_k=\Big\{y\in E: f(y^k,y)\in -C\Big\}.
$$
Finally we define
\begin{equation}\label{ka-9}
x^{k+1}=P_K(w^k).
\end{equation}

Weak convergence of the sequence generated by \eqref{ka-1}--\eqref{ka-9} to a solution of the
vector equilibrium problem was established in \cite{IM-3}. Then the authors in \cite{IM-3},
 performed a minor modification on the above algorithm which ensures strong convergence
of the generated sequence to a solution of VEP$(f,K)$. In Hilbert spaces, this procedure, 
called Halpern's regularization consists of taking a convex combination of a 
given iterate with a fixed point $u\in E$,
where the weight given to $u$ decreases to $0$ with $k$. In Banach spaces, the convex combination must be taken in $E^*$. The strong limit of the generated sequence is the generalized projection of $u$ onto the solution set of the problem.

In this paper, we will consider an extragradient method for solving  vector quasi-equilibrium problems which improves upon \eqref{ka-1}-\eqref{ka-9} in five senses:
\begin{itemize}
\item[i)] We will deal with a rather general class of problems, while \cite{IM-3} only considers vector equilibrium problems.
\item[ii)]  The convergence analysis of the method in \cite{IM-3} requires both  positively weakly upper continuity of $f(\cdot,y)$ for all $y\in E$, and weakly $C$-pseudomonotonicity of $f$, while we avoid using them in this paper.  
\item[iii)] In \cite{IM-3}, the vector valued bifunction $f$ was defined from $E\times E$ to $\re^m$, while we consider a general class of Banach spaces, that is we assume that $f$ is defined  from $E\times E$ to $Y$, where $Y$ is a real Banach space.
\item[iv)] We also show that the boundedness of the sequences generated by our extragradient method implies that the solution set of the vector quasi-equilibrium problem is nonempty, and prove the strong convergence of the generated sequences to a solution of the problem.
\item[v)] In \cite{IM-3}, the authors make intensive use of the auxiliary function $\phi$, but in the current paper we work with a general framework, i.e. the Bregman distance $D_g$ ( see Section 2).      
\end{itemize}  

In this paper, we extend the method in \cite{IM-3} to vector quasi-equilibrium problems, obtaining
an algorithm such that the generated sequence $\{x^k\}$ is strongly convergent to a solution of the problem under minimal assumptions on the bifunction $f$ and the multivalued mapping $T$, and preserving the properties of the method in \cite{IM-3} described in items
(i)-(v) above.

The paper is organized as follows. In Section \ref{s2}, we introduce some preliminary
material  related to the geometry of Banach spaces and vector optimization.
In Section \ref{s3}, we present 
our extragradient method  for solving  
vector quasi-equilibrium problems and prove strong convergence of the generated sequences to a solution of the problem. 
In Section \ref{s4}, we first give some examples of vector quasi-equilibrium problems in several Banach spaces to which  our main theorem can be applied. Then we present some numerical experiments. 

\section{Preliminaries}\label{s2}

Let $E$ be a real Banach space with norm ${\lV\cdot\rV}$. We denote the topological dual of $E$ by $E^*$ and use the notation $\langle x, v \rangle$ for the duality product $v(x)$ of $x\in E$ and $v\in E^*$.
The duality mapping $J:E \to {\cal P}(E^*)$ is defined as
$$
J(x) =\Big \{v \in E^* :\langle x, v \rangle=\lV x\rV^2=\lV v\rV^2\Big\}.
$$
Let $h(x)=\frac{1}{2}\lV x\rV^2$.
It is well known that $h$ is convex and  $J=\partial h$, i.e. $J$ is the subdifferential of
half of the square of the norm.
We assume that  $K\subset E$ is a nonempty, closed and convex set, $Y$ is a
real Banach space containing a closed, convex and pointed cone $C$  with nonempty interior (denoted
as int$(C)$), and
$f:E\times E\rightarrow Y$ is a vector valued bifunction. 

We continue by establishing some standard notation.
The norm, both in $E$ and $Y$, will be denoted as $\lV\cdot\rV$, while the duality coupling between $E$ and
$E^*$, as well as the duality coupling between $Y$ and
$Y^*$ (the topological dual of $Y$), will be denoted as
$\langle\cdot,\cdot\rangle$.
The dual cone $C^+$ of $C$ is defined as $C^+=\{z\in Y^*: \langle y,z\rangle\ge 0, \,\,\forall y\in C\}$.
We define the partial order $\preceq$ in $Y$, induced by the cone $C$, as
$$
y\preceq y'\Longleftrightarrow y'-y \in C,$$
with its associate relation $\prec$, by  $$y\prec y'\Longleftrightarrow y'-y \in {\rm int}(C).$$
We extend $Y$ as $\bar{Y}=Y\cup \{-\infty,+\infty\}$ where a neighbourhood of $+\infty$ is defined as 
a set $N\subset \bar{Y}$ containing $r+C\cup \{+\infty\}$ for some $r\in Y$ and its opposite $-N$ 
is a neighbourhood of $-\infty$. The binary relations $\preceq$ and $\prec$ 
defined in the above are extended to $\bar{Y}$ by 
$$
\forall y\in Y \  \  \ \  \ -\infty \prec y\prec+\infty \  \  \   \   \  -\infty\preceq y\preceq+\infty.
$$
Note that the embedding $Y\subset\bar{Y}$ is continuous and dense. We extend by
continuity every $z\in C^+\setminus\{0\}$ to $\bar{Y}$, by putting $\langle\pm\infty,z\rangle=\pm\infty$. 
Given a set $T\subset \bar{Y}$, we denote its topological closure in the topological space $\bar{Y}$ by $\bar{T}$.
To a given set $T\subset\bar{Y}$, we associate the following set:
$$
{\rm inf}^C_w(T)=\Big\{y\in \bar{T} | \not\exists z\in T: z\prec y\Big\}.
$$
Given $S\subset E$ and $G: S\to Y\cup \{+\infty\}$, the point $a\in E$ is called 
{\it weakly efficient} if $a\in S$ and $G(a)\in {\rm inf}^C_w(G(S))$.
We denote as ${\rm argmin}^C_w\{G(x)| x\in S\}$ the set of weakly efficient points. We observe that 
$$
{\rm argmin}^C_w\Big\{G(x)| x\in S\Big\}=S\cap G^{-1}({\rm inf}^C_w(G(S))).
$$
\begin{definition}
A map $G:E\rightarrow Y\cup \{+\infty\}$ is called C-convex whenever 
$$
G(tx+(1-t)y)\preceq tG(x)+(1-t)G(y), \ \ \ \ \ \forall x,y \in E \ \ \ {\rm and}\ \ \ \forall  t\in [0, 1].
$$
\end{definition}


A Banach space $E$ is said to be {\it strictly convex} if $\lV\frac{x + y}{2}\rV<1$ for all $x, y\in E$
with $\lV x\rV =\lV y\rV = 1$ and $x \neq y$. It is said to be {\it uniformly convex} if for each $\varepsilon\in (0,2]$, 
there exists $\delta>0$ such that for all $x, y\in E$
with $\lV x\rV =\lV y\rV = 1$ and $\lV x- y\rV\ge\varepsilon$,  it holds that
$\lV\frac{x + y}{2}\rV<1-\delta$.
It is known that  uniformly convex Banach spaces are
reflexive and strictly convex. 

A Banach space $E$ is said
to be {\it smooth} if 
\begin{equation}\label{smooth}
\lim_{t\to 0} \frac{\lV x+ty\rV -\lV x\rV}{t}
\end{equation}
exists for all $x, y \in B_1(0) = \{z \in E : \lV z\rV = 1\}$. It is said to be {\it uniformly smooth} if the limit in \eqref{smooth} is
attained uniformly for $x, y \in B_1(0)$. It is well known that the spaces $L^p$ $(1<p<+\infty)$ and the Sobolev spaces $W^{k,p}$ ($1<p<+\infty$) 
are both uniformly convex and uniformly smooth. 

Now we recall some properties of Bregman distance which will be used in this paper (see \cite{BuI}, \cite{BIR} and \cite{IN}). We consider an auxiliary function 
$g : E \to \re$, which is strictly convex, lower semicontinuous, and G\^ateaux differentiable. We will denote the family of such functions by $\cal{F}$. The G\^ateaux derivative of $g$ will be denoted by $g'$.

\begin{definition}
Let $g : E \to \re$ be a convex and G\^ateaux differentiable function.\\
i) The Bregman distance with respect to $g$ is the function $D_g : E\times E \to \re$, defined by
\begin{equation}\label{breg}
D_g (x, y) = g(x)-g(y)-\langle x-y,g'(y)\rangle.
\end{equation}
ii) The modulus of total convexity of $g$ is the function 
$v_g : E\times [0,+\infty)\to [0,+\infty)$
defined by $v_g(x,t)=\inf\Big\{D_g(y, x):y \in E, \ \|y-x\| = t\Big\}$.\\
iii) $g$ is said to be a totally convex function at $x\in E$ if $v_g(x,t)>0$ for all $t>0$.\\
iv) $g$ is said to be a totally convex function if $v_g(x,t)>0$ for all $t>0$ and all $x\in E$.\\
v) $g$ is said to be a uniformly totally convex function on $B\subset E$ if $\inf_{x\in A}v_g(x,t)>0$ for all $t>0$ and all bounded subsets $A\subset B$.
\end{definition}
It is worthwhile mentioning that $D_g (x, y) =\|x-y\|^2$
 whenever $g(x)=\|x\|^2$ and $E$ is a Hilbert space.  
Now we assume some additional conditions on $g\in \cal{F}$, which will be needed for the convergence analysis of our algorithm.\\
\\
H1: The level sets of $D_g (x, \cdot)$ are bounded for all $x\in E$.\\
H2: $g$ is uniformly totally convex on $E$.\\
H3: $g'$ is uniformly continuous on bounded subsets of $E$.\\
H4: $\lim_{\|x\|\to \infty}\Big(g(x)-\rho\|x-z\|\Big)=\infty$ for all fixed $z\in E$ and $\rho>0$.

\begin{proposition}\label{lph5}{\rm (\cite{IN}, Proposition 2.3)}
If $E$ is a uniformly smooth and uniformly convex Banach space, then
$g(x)=r\|x\|^s$ satisfies H1--H4 for all $r>0$ and all $s>1$.
\end{proposition}

It is well known that when $E$ is smooth, the duality operator $J$ is single valued.
Let $E$ be a smooth Banach space. We define $\phi:E\times E\to\re$ by
\begin{equation}\label{phi-function}
\phi(x,y)=\lV x\rV^2-2\langle x, J(y) \rangle +\lV y\rV^2.
\end{equation}
This function can be seen as a ``distance-like" function, better conditioned than the square of the
metric distance, namely $\lV x-y\rV^2$; see  e.g. \cite{Alb}, \cite{KaT} and \cite{Rei}.
In Hilbert spaces, where the duality mapping $J$ is the identity operator, it holds that $\phi(x,y)=\lV x-y\rV^2$. Moreover, if we define  $g(x)=\lV x\rV^2$ for all $x\in E$, then $D_g(x,y)=\phi(x,y)$. In the sequel, we will need the following properties of the Bregman distance $D_g$.
\begin{proposition}\label{zero}{\rm (\cite{IG}, Proposition 5)}
Suppose that $g\in\cal{F}$ satisfies H2.  Let $\{x^k\}$ and $\{y^k\}$ be  two sequences in $E$.  
If  $\lim_{k\to\infty}D_g(x^k,y^k)= 0$ and either $\{x^k\}$ or $\{y^k\}$ is bounded, then 
$\lim_{k\to\infty}\lV x^k-y^k\rV$ $=0$. 
\end{proposition}

\begin{proposition}\label{lem-proj}{\rm (\cite{BuI}, Page 70)}
Let $K\subset E$ be nonempty, closed and convex, and $g\in\cal{F}$ be a totally convex function on
$E$ satisfying H1--H2. Consider $x\in E$, then there exists a unique $\bar{x}\in K$ such that 
$$\bar{x}={\rm argmin}_{y\in K}D_g(y,x).$$
We denote $\bar{x}=\Pi^g_K(x)$ and call $\Pi^g_K$ the Bregman projection operator from $E$ onto $K$.
Moreover, $\bar{x}=\Pi^g_K(x)$ if and only if 
$$
\langle z-\bar{x},g'(x)-g'(\bar{x})\rangle\le 0
$$  
for all $z\in K$.
\end{proposition}
\begin{proposition}\label{uni-bound}{\rm (\cite{IG}, Proposition 4)}
If $g$ satisfies H3, then both $g$ and $g'$ are bounded on bounded subsets of $E$.
\end{proposition}

Now we introduce some notations and definitions that will be used in the sequel. For a sequence 
$\{x^k\}$ in $E$, we denote strong convergence of $\{x^k\}$ to $x\in E$ by $x^k\rightarrow x$, 
and weak convergence by $x^k\rightharpoonup x$. In the following definitions, we assume that  $K\subset E$ is a nonempty, closed and convex set, and $g\in {\cal F}$
is a totally convex function on $E$ satisfying H1--H2. 

\begin{definition}
 We say that $T:K\rightarrow K$ is a quasi $D_g$-nonexpansive mapping whenever ${\rm Fix}(T)\not=\emptyset$ and $D_g(p,Tx)\leq D_g(p,x)$ for all $(p,x)\in {\rm Fix}(T)\times K$.
\end{definition}
\begin{definition}
Let $T(\cdot)$ be a multivalued mapping from $K$ into ${\cal P}(K)$ such that 
for all $x\in K$, $T(x)$  is a nonempty, closed and convex subset of $K$. We say that $T(\cdot)$ is quasi $D_g$-nonexpansive  whenever the mapping $S(\cdot)=\Pi^g_{T(\cdot)}(\cdot)$ is quasi $D_g$-nonexpansive where $\Pi^g$ is the Bregman projection.
\end{definition}

\begin{definition}\label{defdemic}
 The multivalued mapping $T(\cdot)$ from $K$ into ${\cal P}(K)$ is said to be demiclosed, if whenever   $x^k\rightharpoonup \bar{x}$ and $\lim_{k\to\infty}d(x^k,T(x^k))=0$, then $\bar{x}\in{\rm Fix}(T)$.
 \end{definition}

\begin{proposition}\label{FTCC}
 If $T(\cdot):K\to {\cal P}(K)$ is a  quasi $D_g$-nonexpansive  mapping, then ${\rm Fix}(T)$ is closed and convex.
\end{proposition}
\begin{proof}
 Let $p_1,p_2\in {\rm Fix}(T)$ and define $p_t=tp_1+(1-t)p_2$ where $t\in [0,1]$. In order to prove the convexity of ${\rm Fix}(T)$, we must show that $p_t\in T(p_t)$.  Let $Sp_t:=\Pi^g_{T(p_t)}(p_t)$ where $\Pi^g$ is the Bregman projection. Note that $D_g$ is nonnegative, then by the definition of the Bregman distance, we have 
\begin{align*}\label{ft-}
0&\leq D_g(p_t,Sp_t)=g( p_t)-g(Sp_t)-\langle p_t-Sp_t,g'(Sp_t)\rangle\\
&=g( p_t)-tg( p_1)-(1-t)g( p_2)+t(g( p_1)-g(Sp_t)-\langle p_1,g'(Sp_t)\rangle+\langle Sp_t,g'(Sp_t)\rangle)\\
&+(1-t)(g( p_2)-g(Sp_t)-\langle p_2,g'(Sp_t)\rangle+\langle Sp_t,g'(Sp_t)\rangle)\\
&=g( p_t)-tg( p_1)-(1-t)g( p_2)+tD_g(p_1,Sp_t)+(1-t)D_g(p_2,Sp_t)\\
&\leq g( p_t)-tg( p_1)-(1-t)g( p_2)+tD_g(p_1,p_t)+(1-t)D_g(p_2,p_t)\\
&=D_g(p_t,p_t)=0.
\end{align*}
Therefore $D_g(p_t,Sp_t)=0$. Now Proposition \ref{zero} shows that $ Sp_t=p_t$. Since $Sp_t=\Pi^g_{T(p_t)}(p_t)$,  $p_t\in T(p_t)$, i.e. ${\rm Fix}(T)$ is convex.\\
Now we show that ${\rm Fix}(T)$ is closed. Let  $\{p^k\}\subset {\rm Fix}(T)$ be such that $p^k\to p$, and let $Sp=\Pi^g_{T(p)}(p)$. 
 Since $D_g(p^k,Sp)\leq D_g(p^k,p)$ for all $k$, we have $\lim_{k\to \infty}D_g(p^k,Sp)=D_g(p,Sp)=0$. Then Proposition \ref{zero} implies that   $p\in {\rm Fix}(S)$. Therefore
$p\in {\rm Fix}(T)$, i.e. ${\rm Fix}(T)$ is closed.
\end{proof}

\begin{definition}
 The multivalued mapping $T(\cdot)$ from $K$ to itself is called lower semicontinuous at each $\bar{x}\in K$,  whenever we have $\{x^k\}\subset K$ and $x^k\rightarrow \bar{x}$, then for any $\bar{y}\in T(\bar{x})$, there is a sequence $\{y^k\}$ with $y^k\in T(x^k)$ for all k, such that 
$y^k\rightarrow \bar{y}$ as $k\to \infty$.
\end{definition}
In the following, we give an example of a multivalued mapping which is quasi $D_g$-nonexpansive, demiclosed and  lower semicontinuous at each $\bar{x}\in K$.
\example\label{fir-ex-1}
Define  $T(\cdot):K\to {\cal P}(K)$ as $T(x)=B(0, \|x\|)$ where $B(0,\|x\|)$ denotes the closed ball of radius $\|x\|$ centered at $0$. It is easy to see that $T$ is demiclosed and quasi $D_g$-nonexpansive with $g(\cdot)=\|\cdot\|^2$. Now we show that $T$ is lower semicontinuous at each $\bar{x}\in K$. Suppose that $x^k\rightarrow \bar{x}$ and $\bar{y}\in T(\bar{x})$. Then if $\bar{x}=0$, we have $\bar{y}=0$ and hence we define $y^k=0$ for all $k$, also we have $y^k\in T(x^k)$. In the sequel, if $\bar{x}\not=0$ and $\bar{y}\in T(\bar{x})$, we define $y^k=\frac{\langle x^k, J\bar{x}\rangle}{\|\bar{x}\|^2}\bar{y}$. It easy to see that $y^k\in T(x^k)$. Therefore in both cases we have $y^k\in T(x^k)$ such that
$y^k\to \bar{y}$.

Now we introduce some assumptions on the vector valued bifunction $f:E\times E\rightarrow Y$  and the multivalued mapping $T$, that we will need for the convergence analysis.

\begin{itemize}
\item[B1:] $f(x,x)=0$ for all $x\in E$,
\item[B2:] $f(\cdot,\cdot):E\times E\to Y$ is uniformly continuous on bounded sets,
\item[B3:] $f(x,\cdot):E\to Y$ is $C$-convex for all $x\in E$.
\item[B4:] $T(\cdot):K\to {\cal P}(K)$ is a multivalued mapping with nonempty, closed and convex values, demiclosed, lower semicontinuous and quasi $D_g$-nonexpansive, where $g\in {\cal F}$.
\end{itemize}

 We also mention that for the sequences generated by our algorithm in Section \ref{s3} to be well defined and bounded,  we will assume that 
$$DS(f,T):=\Big\{x\in T(x) : f(y,x)\in -C, \ \forall y \in K \Big\}\not=\emptyset.$$
However, if the sequences generated by the algorithm SEML, introduced in Section \ref{s3}, are well defined and bounded, then we show that the vector quasi-equilibrium problem has a solution.

A vector valued function $G:E\rightarrow Y\cup \{+\infty\}$ is called positively lower semicontinuous, 
if for every $z\in C^+$ the extended scaler function $x\mapsto\langle G(x), z\rangle$ is 
lower semicontinuous. Also we say that
$G$ is positively upper semicontinuous whenever $-G$ is positively lower semicontinuous.

Now we recall an essential theorem from \cite{BIS}, which is needed in the next sections.
\begin{theorem}\label{theo-21}
If $S\subset E$ is a convex set and $G:S \to Y\cup \{+\infty\}$ is a $C$-convex proper map, then 
$$
{\rm argmin}^C_w \Big\{G(x)\mid x\in S\Big\}=\bigcup_{z\in C^+\setminus\{0\}}{\rm argmin}\Big\{\langle G(x),z\rangle\mid x\in S\Big\}.$$
\end{theorem}

  We also need the following result from \cite{phelps}.
\begin{proposition}{ \rm \cite{phelps}}\label{theo-phelps}
 Suppose that $f$ and $g$ are  proper, convex and lower semicontinuous functions on the Banach space $E$ and that there is a point in $D(f)\cap D(g)$ where one of them is continuous. Then
 $$\partial(f + g)(x) = \partial f(x)+\partial g(x), \ \ \ \ \ \  x\in D(\partial f)\cap D(\partial g).$$
\end{proposition}

We recall now some properties of the solution set of dual vector quasi-equilibrium problems.

\begin{proposition}\label{pp1}
Assume that $T(\cdot):K\to {\cal P}(K)$ is a multivalued mapping and $f:E\times E\to Y$ satisfies B1, that $f(\cdot,y)$ is positively upper semicontinuous for all $y\in E$ 
and that $f(x,\cdot)$ is $C$-convex for all $x\in E$. Then $DS(f,T)\subset S(f,T)$. 
\end{proposition}
\begin{proof}
Take $x^*\in DS(f,T)$ and define $p_t=tx^*+(1-t)y$ with $t\in (0,1)$  
and $y\in T(x^*)$. Take any $c\in C^+\setminus\{0\}$. B1 and $C$-convexity of $f(p_t,\cdot)$ imply that 
\begin{equation}\label{dss}
0=\langle f(p_t,p_t), c\rangle\leq t\langle f(p_t,x^*), c\rangle+(1-t)\langle f(p_t,y), c\rangle.
\end{equation}
Since $\langle f(p_t,x^*), c\rangle\leq 0$, \eqref{dss} implies that 
\begin{equation}\label{ea1}
\langle f(p_t,y), c\rangle\geq 0.
\end{equation} 
Since $\langle f(\cdot,y), c\rangle$ is upper semicontinuous, taking limsup with $t\to 1$ in 
\eqref{ea1} gives  $\langle f(x^*,y), c\rangle\geq 0$. Hence $f(x^*,y)\not \in -{\rm int}(C)$. 
Since $y\in T(x^*)$ is arbitrary, we get $DS(f,T)\subset S(f,T)$.
\end{proof}

\begin{corollary}\label{ca1}
Under B1--B3, $DS(f,T)\subset S(f,T)$.
\end{corollary}
\begin{proof}
Elementary.
\end{proof}

\begin{proposition}\label{pp2}
If  $T(\cdot):K\to {\cal P}(K)$ is a quasi $D_g$-nonexpansive mapping, and $f(x,\cdot)$ is $C$-convex and positively lower semicontinuous for all $x\in E$, then $DS(f,T)$ is closed and convex.
\end{proposition}
\begin{proof}
Note that ${\rm Fix}(T)$ is closed and convex by Proposition \ref{FTCC}. Now, take $\bar{x} ,x^*\in DS(f,T)$ and define $x_t=tx^*+(1-t)\bar{x}$ with $t\in (0,1)$. Take any $c\in C^+\setminus\{0\}$. 
By  $C$-convexity of $f(x,\cdot)$, we have  
\begin{equation}\label{dscon}
\langle f(x,x_t), c\rangle\leq t\langle f(x,x^*), c\rangle+(1-t)\langle f(x,\bar{x}), c\rangle\leq0, 
\end{equation}
for all $x\in K$.
Since $c\in C^+\setminus\{0\}$ is arbitrary, it follows that $x_t\in DS(f,T)$, i.e. $DS(f,T)$ is convex.
Closedness of 
$DS(f,T)$  
follows from positive lower semicontinuity of $f(x,\cdot)$ for all $x\in E$. 
\end{proof}

\begin{corollary}\label{ca2}
Under B1--B4, $DS(f,T)$ is closed and convex.
\end{corollary}
\begin{proof}
Follows from Propositions \ref{FTCC} and \ref{pp2}.
\end{proof}

\section{Extragradient method with linesearch and strong convergence}\label{s3}

\noindent In this section, we study the strong convergence of the
sequence generated by a {\bf Strongly convergent variant of the Extragradient Method with Linesearch (SEML)} to approximate a solution of the vector quasi-equilibrium problem. We propose a regularization procedure on the extragradient method which ensures the strong convergence
of the generated sequence to a solution of the problem. 
We will assume in the sequel that $E$ is a Banach space and $K\subset E$ is nonempty closed and convex,  a real Banach space  $Y$ containing a closed, convex and pointed cone $C$  with nonempty interior, and  that $f:E\times E\to Y$ is a vector valued bifunction, $T(\cdot):K\to {\cal P}(K)$ is a multivalued  mapping, $g\in {\cal F}$ is a  function on $E$ satisfying H1--H4,  and the assumptions B1--B4 are satisfied. For the sake of definiteness and boundedness  of
the iterative sequences $\{x^k\}$, $\{v^k\}$ and $\{w^k\}$ generated by the following algorithm, we assume that $DS(f,T)\neq\emptyset$. However, we will show later that if the sequences generated by the algorithm are bounded, then $S(f,T)\neq\emptyset$.
First we give the formal definition of Algorithm SEML.

\bigskip

\noindent{\bf 1. Initialization}:

Fix $v^0\in K$ and  $\delta, \theta\in(0,1)$. 
Take $\hat\beta$, $\tilde\beta$ satisfying
$0<\hat{\beta}\leq\tilde{\beta}$, and consider a sequence
$\{\beta_k\}\subset[\hat\beta,\tilde\beta]$ and a sequence 
$\{\gamma_k\}\subset [\varepsilon, 1]$ for some $\varepsilon \in (0, 1]$.
Also, take a sequence $\{e^k\}\subset {\rm int}(C)$ such that $e^k\to \bar{e} \in {\rm int}(C)$.

\noindent{\bf 2. Iterative step}:  Given $v^k$, define
\begin{equation}\label{ka-1s}
 x^k=\Pi^g_{T(v^k)}(v^k).
\end{equation}
\noindent{\bf 2. Iterative step}:
\begin{equation}\label{ka-2s}
 z^k \in {\rm argmin}^C_w\left\{\beta_kf(x^k,y)+
g(y)e^k-\langle y, g'(x^k) \rangle e^k : y\in T(v^k) \right\}.
\end{equation}
If $z^k=v^k$ stop. Otherwise, let
\begin{equation}\label{ka-4s}
\ell(k)=\min\left\{\ell \geq 0: -\beta_kf(y^{\ell},x^k)+\beta_kf(y^{\ell},z^k)+ 
\delta D_g(z^k,x^k)e^k\not\in {\rm int}(C) \right\},
\end{equation}
with
\begin{equation}\label{ka-5s}
y^{\ell}=\theta^{\ell}z^k+(1-\theta^{\ell})x^k.
\end{equation}
Set
\begin{equation}\label{ka-6s}
\alpha_k:=\theta^{\ell(k)},
\end{equation}
\begin{equation}\label{ka-7s}
y^k:=y^{\ell(k)}=\alpha_k z^k+(1-\alpha_k)x^k.
\end{equation}
Define  
\begin{equation}\label{ka-8sh}
H_k:=\Big\{y\in E: f(y^k,y)\in -C\Big\},
\end{equation}
If $k=0$, set $K_0=K\cap H_0$. Otherwise, let 
\begin{equation}\label{ka-8sc}
K_k=K_{k-1}\cap H_k.
\end{equation}
\begin{equation}\label{ka-8s}
w^k=\Pi^g_{K_k}(x^k).
\end{equation}
Determine the next approximation $v^{k+1}$ as
\begin{equation}\label{ka-11s}
v^{k+1} = \Pi^g_{L_k\cap M_k\cap N_k}(v^0),
\end{equation}
where
\begin{equation}\label{ka-12s}
L_k=\Big\{ z\in E: \langle z-x^k,g'(x^k)-g'(w^k)\rangle \leq -\gamma_k D_g(x^k,w^k)\Big\},
\end{equation}
\begin{equation}\label{ka-13s}
M_k=\Big\{ z\in E: \langle z-v^k,g'(v^k)-g'(x^k)\rangle \leq -\gamma_k D_g(v^k,x^k)\Big\},
\end{equation}
\begin{equation}\label{ka-14s}
N_k= \Big\{ z\in E: \langle z-v^k, g'(v^0)-g'(v^k)\rangle \leq0\Big\}.
\end{equation}
\bigskip

We proceed now to the convergence analysis of Algorithm SEML. The proof of the main theorem is divided into several Lemmas and Propositions. In order to establish the strong convergence of the sequences $\{x^k\}$ and $\{v^k\}$  generated by Algorithm SEML, both to a solution of the problem,  we need some intermediate results. 

\begin{theorem} \label{main-strong}
Assume that $f$ is a vector valued bifunction, $g\in {\cal F}$ is a  function on $E$ satisfying H1--H4, $T(\cdot)$ is a multivalued  mapping from $K$ to ${\cal P}(K)$ and the assumptions B1--B4 are satisfied.\\
 i) If $DS(f,T)\neq\emptyset$, then the sequences $\{x^k\}$, $\{v^k\}$ and $\{w^k\}$ generated by Algorithm SEML are well defined and bounded.\\
 ii) If the sequences generated by the algorithm are well defined and bounded, then the sequences $\{v^k\}$ and $\{x^k\}$, both converge strongly to an element of $S(f,T)$, which is therefore nonempty.
\end{theorem}

We will give the proof of Theorem \ref{main-strong} at the end of this section, after proving the intermediary steps needed for the proof.

\begin{proposition}\label{well-1s}
The sequence $\{ z^k\}$ generated by Algorithm SEML is well defined.
\end{proposition}
\begin{proof}
Take any $c\in C^+\setminus \{0\}$. Since $e^k \in {\rm int}(C)$, it follows from the 
definition of $C^+$ that $\langle e^k,c\rangle>0$. Define $\psi :E\to \re \cup\{+\infty\}$ as 
\begin{equation}\label{well-2}
\psi(y)= \beta_k \langle f(x^k,y),c\rangle+
g(y) \langle e^k,c\rangle-\langle y, g'(x^k) \rangle \langle e^k,c\rangle.
\end{equation}
It is easy to see that $\psi$ is proper, convex and lower semicontinuous. The subdifferential of $\psi$ is maximal monotone, 
and hence onto by Corollary 3.7 of \cite{phelps}. Thus $\partial\psi$  has some zero, which is a minimizer of $\psi$. 
In view of Theorem \ref{theo-21} such minimizer satisfies \eqref{ka-2s} and can be taken as $z^k$.
\end{proof}

\begin{proposition}\label{pp}
Assume that $f$ satisfies B1--B3. Take $v\in K$, $x\in T(v)$, $\beta\in\re^+$ and $e\in {\rm int}(C)$. If  
\begin{equation}\label{KOR1} 
z\in {\rm argmin}^C_w\Big\{\beta f(x,y)+
g(y)e-\langle y, g'(x) \rangle e : y\in T(v)\Big\} 
\end{equation}
then there exists $c\in C^+\setminus\{0\}$ such that 
$$\langle y-z,g'(x)-g'(z)\rangle \langle e, c\rangle \leq
\beta[\langle f(x,y), c\rangle-\langle f(x,z), c\rangle], \ \ \ \ \ \ \forall y\in T(v).$$
\end{proposition}
\begin{proof}
Let $N_{T(v)}(z)$ be the normal cone of $T(v)$ at $z\in {T(v)}$, i.e. 
$N_{T(v)}(z)=\{v^*\in E^* : \langle y-z, v^*\rangle\leq 0, \forall y\in T(v)\}$.
Since $z$ solves the vector optimization problem in \eqref{KOR1}, in view of Theorem \ref{theo-21} there exists 
$c\in C^+\setminus \{0\}$ such that $z$ satisfies the
first order optimality condition, given by 
$$
0\in \partial\Big\{\langle \beta f(x,\cdot), c\rangle+
g(\cdot)\langle e, c\rangle-\langle \cdot, g'(x)\rangle \langle e, c\rangle\Big\}(z)+N_{T(v)}(z).
$$
Thus, in view of the definition of $g'$, and by Proposition \ref{theo-phelps}, there exist 
$w\in \partial \langle f(x,\cdot), c\rangle(z)$ 
and 
$\bar{w}\in N_{T(v)}(z)$ 
such that 
$$
0=\beta w+\langle e, c\rangle g'(z)-\langle e, c\rangle g'(x)+\bar{w}.
$$
Therefore, since   
$\bar{w}\in N_{T(v)}(z)$, we have 
$\langle y-z,-\beta w-\langle e, c\rangle g'(z)+\langle e, c\rangle g'(x)\rangle\le 0$, so that,
using the fact that
$w\in \partial \langle  f(x,\cdot), c\rangle(z)$, we get
\begin{equation}\label{eq-0}
\langle y-z,g'(x)-g'(z)\rangle\langle e, c\rangle\le\beta\langle y-z,w\rangle\le
\beta\langle  f(x,y), c\rangle-\beta\langle  f(x,z), c\rangle.
\end{equation}

\end{proof}

\begin{corollary}\label{c1} Assume that  $\{x^k\}$ and $\{z^k\}$ are the sequences generated by Algorithm SEML. Then there 
exists $\{c^k\}\subset C^+\setminus\{0\}$ such that
$$
\langle y-z^k,g'(x^k)-g'(z^k)\rangle \langle e^k, c^k\rangle \le
\beta_k\left[\langle f(x^k,y), c^k\rangle-\langle f(x^k,z^k), c^k\rangle\right] \ \ \ \ \ \ \forall y\in T(v^k).
$$
\end{corollary}
\begin{proof} Follows from Proposition \ref{pp} and \eqref{ka-2s}.
\end{proof}

\begin{proposition}\label{stop-1}
If Algorithm SEML stops at the $k$-th iteration, then $x^k$ is a solution of VQEP($f,T$).
\end{proposition}
\begin{proof}
If $z^k=v^k$, since $z^k\in T(v^k)$, we get $z^k=x^k$ by \eqref{ka-1s}. Hence Corollary \ref{c1} implies that $\langle f(x^k,y), c^k\rangle \ge 0$ for all $y\in T(x^k)$. 
Since $c^k\in C^+\setminus\{0\}$, we have $f(x^k,y)\not\in -{\rm int}(C)$ for all $y\in {T(x^k)}$.
\end{proof}

\begin{proposition}\label{lk-1}
The following statements hold for Algorithm SEML.\\
i) $\ell(k)$ is well defined, (i.e. the Armijo-type search for $\alpha_k$ is finite), and consequently the 
same holds for the sequence $\{y^k\}$.\\
ii)  If $x^k\not=z^k$, then $f(y^k,x^k)\not\in -C$.
\end{proposition}

\begin{proof}
i) We proceed inductively, i.e. we assume that $v^k$ is well defined, and proceed to
establish that the same holds for $v^{k+1}$. Note that $z^k$ is well defined by Proposition \ref{well-1s} and also $x^k$ is well defined by Proposition \ref{lem-proj}. It suffices to
check that $\ell(k)$ is well defined. Assume by contradiction that
\begin{equation}\label{}
 -\beta_kf(y^{\ell},x^k)+\beta_kf(y^{\ell},z^k)+ \delta D_g(z^k,x^k)e^k\in {\rm int}(C)
\end{equation}
for all $\ell$. Since $c^k\in C^+\setminus\{0\}$,  we have
\begin{equation}\label{www}
\beta_k[\langle f(y^{\ell},x^k),c^k\rangle-\langle f(y^{\ell},z^k),c^k\rangle]
< \delta D_g(z^k,x^k)\langle e^k,c^k\rangle
\end{equation}
for all $\ell$. Note that the sequence $\{y^\ell\}$ is strongly convergent to $x^k$. In view of B2, taking limits in \eqref{www}
as $\ell\to +\infty$, 
\begin{equation}\label{w1}
\beta_k[\langle f(x^k,x^k),c^k\rangle-\langle f(x^k,z^k),c^k\rangle]
\le \delta D_g(z^k,x^k)\langle e^k,c^k\rangle.
\end{equation}
Since $x^k\in T(v^k)$ by \eqref{ka-1s}, we apply Corollary \ref{c1} with $y=x^k$ in \eqref{w1}, obtaining 
\begin{equation}\label{w2}
\langle x^k-z^k,g'(x^k)-g'(z^k)\rangle\leq\delta D_g(z^k,x^k).
\end{equation}
In view of the definition of $ D_g$, \eqref{w2} implies that 
\begin{equation}
 D_g(z^k,x^k)+ D_g(x^k,z^k)\le\delta D_g(z^k,x^k).
\end{equation}
Since $\delta\in(0,1)$, we get
$ D_g(x^k,z^k)< 0$,
contradicting the nonnegativity of $ D_g$.\\
ii)  Assume that $f(y^k,x^k)\in -C$. Note that, using B1, B3 and \eqref{ka-7s}, we have
$$
0=f(y^k,y^k)\preceq \alpha_kf(y^k,z^k)+(1-\alpha_k)f(y^k,x^k).
$$
Since $-(1-\alpha_k)f(y^k,x^k)$ and $\alpha_kf(y^k,z^k)+(1-\alpha_k)f(y^k,x^k)$
belong to $C$, and $C$ is a convex cone, 
we conclude that $f(y^k,z^k)\in C$. Therefore 
\begin{equation}\label{e500}
-\beta_kf(y^k,x^k)+\beta_kf(y^k,z^k)+ \delta D_g(z^k,x^k)e^k\in {\rm int}(C),
\end{equation}
which contradicts \eqref{ka-4s}--\eqref{ka-7s}. Note that the inclusion  
in \eqref{e500} is due to the fact that 
$x^k\neq z^k$ and $e^k\in {\rm int}(C)$. 
\end{proof}

In order to prove the strong convergence of the sequences $\{x^k\}$ and $\{v^k\}$ generated by the algorithm, we need the following lemmas.

\begin{lemma}\label{xn}
If $DS(f,T)\not=\emptyset$, then $DS(f,T)\subset L_k \cap M_k\cap N_k$. Therefore the sequences $\{v^k\}$, $\{w^k\}$ and $\{x^k\}$ are well defined.
\end{lemma}

\begin{proof} The proof is by induction. Note that $DS(f,T)$, $L_k$, $M_k$ and $N_k$ are closed and
convex. We first show that $DS(f,T)\subset L_k \cap M_k\cap N_k$ for all $k\geq0$. Putting $$D_k=\Big\{ z\in E:  D_g(z,w^k)\leq  D_g(z,x^k)\Big\}=\Big\{ z\in E:  \langle z-x^k,g'(x^k)-g'(w^k)\rangle\leq - D_g(x^k,w^k) \Big\}$$ and
$$F_k=\Big\{ z\in E:  D_g(z,x^k)\leq  D_g(z,v^k)\Big\}=\Big\{ z\in E:  \langle z-v^k,g'(v^k)-g'(x^k)\rangle\leq - D_g(v^k,x^k) \Big\}.$$
By $\gamma_k\in [\varepsilon, 1]$, we get $D_k\subset L_k$ and $F_k\subset M_k$. Let $x^*\in DS(f,T)$, note that $x^*\in H_{k}$ for all $k$, we also have $w^{k}=\Pi^g_{K_k}(x^k)$ by \eqref{ka-8s}. Now Proposition \ref{lem-proj} implies that 
$$
\langle x^*-w^{k},g'(x^k)-g'(w^k)\rangle \leq 0,
$$
or equivalently,
\begin{equation}\label{eq1h}
 D_g(w^{k},x^{k})+ D_g(x^*,w^{k})- D_g(x^*,x^{k})\le 0.
\end{equation}
Therefore we have
\begin{equation}\label{w-bound}
 D_g(x^*,w^k)\leq  D_g(x^*,x^k),
\end{equation}
which implies that $DS(f,T)\subset D_k$ for all $k\geq 0$.

On the other hand, since $x^{k}=\Pi^g_{T(v^k)}(v^k)$ and $\Pi^g_{T(\cdot)}(\cdot)$ is a quasi $ D_g$-nonexpansive mapping, we have
\begin{equation}\label{v-bound}
  D_g(x^*,x^k)\leq  D_g(x^*,v^k)
\end{equation}
for all $x^*\in DS(f,T)$. Therefore $DS(f,T)\subset D_k\cap F_k$ for all $k\geq 0$, that implies $DS(f,T)\subset L_k\cap M_k$ for all $k\geq0$. Next, we show that $DS(f,T)\subset L_k \cap M_k\cap N_k$, for all $k\geq 0$, by the induction. Indeed, we have $DS(f,T)\subset L_0 \cap M_0\cap N_0$, because $N_0=E$. Assume that $DS(f,T)\subset L_k \cap M_k\cap N_k$ for some $k\geq 0$. Since
$v^{k+1} = \Pi^g_{L_k \cap M_k\cap N_k}(v^0)$, we have by proposition \ref{lem-proj} that
$$\langle z-v^{k+1}, g'(v^0)-g'(v^{k+1})\rangle \leq0, \  \  \forall z\in L_k \cap M_k\cap N_k.$$
Since $DS(f,T)\subset L_k \cap M_k\cap N_k$, we have
$$\langle z-v^{k+1}, g'(v^0)-g'(v^{k+1})\rangle \leq0, \  \  \forall z\in DS(f,T).$$
Now, since  $\langle z-v^{k+1}, g'(v^0)-g'(v^{k+1})\rangle\leq0,\ \  \forall z \in DS(f,T)$, the definition of $N_{k+1}$ implies that $DS(f,T)\subset N_{k+1}$, and so
$DS(f,T)\subset L_k \cap M_k\cap N_k$ for all $k\geq 0$. Finally, since $DS(f,T)$ is nonempty, we get $L_k \cap M_k\cap N_k$ is nonempty,  therefore $v^{k+1}$ is well defined. Now, it is clear that the sequences $\{x^k\}$ and $\{w^k\}$ are well defined.
\end{proof}
\begin{lemma}\label{bound-vk}
If $DS(f,T)\not=\emptyset$, then the sequences $\{x^k\}$, $\{v^k\}$ and $\{w^k\}$ generated by Algorithm SEML are bounded.
\end{lemma}
\begin{proof}
From the definition of $N_k$, we have $v^k=\Pi^g_{N_k}(v^0)$. Let $x^*\in DS(f,T)$.  Since $DS(f,T) \subset N_k$ by Lemma \ref{xn}, and $\Pi^g_{N_k}$ is the Bregman projection onto $N_k$, we have $\langle x^*- v^k, g'(v^0)-g'(v^k)\rangle \leq0$ by Proposition \ref{lem-proj}, which implies 
\begin{equation}\label{}
D_g(x^*,v^k)\leq D_g(x^*,v^0).
\end{equation} 
Thus, the sequence $\{v^k\}$ is bounded by H1.\\
Also, since $x^{k}=\Pi^g_{T(v^k)}(v^k)$ and $\Pi^g_{T(\cdot)}(\cdot)$ is a quasi $D_g$-nonexpansive mapping, we have
\begin{equation}\label{}
 D_g(x^*,x^k)\leq D_g(x^*,v^k).
\end{equation}
 Therefore the boundedness of  the sequence $\{v^k\}$
implies that the sequence $\{x^k\}$ is bounded by H1.\\
In the sequel, since $x^*\in DS(f,T)$, we have $f(y^k, x^*)\in -C$, hence \eqref{ka-8sh} shows that  $x^*\in H_k$ for all $k$, therefore $x^*\in K_k$ for all $k$.
On the other hand, since  $w^{k}=\Pi^g_{K_k}(x^k)$ by \eqref{ka-8s}, Proposition \ref{lem-proj} implies that 
$$
\langle x^*-w^k,g'(x^k)-g'(w^k)\rangle \leq 0.
$$
Hence
\begin{equation}\label{}
D_g(w^k,x^k)+D_g(x^*,w^k)-D_g(x^*,x^k)\le 0.
\end{equation}
Therefore we have
\begin{equation}\label{}
D_g(x^*,w^k)\leq D_g(x^*,x^k).
\end{equation}
Since the sequence $\{x^k\}$ is bounded, it follows from H1 that the sequence $\{w^k\}$ is bounded too.\\

\end{proof}

\begin{lemma}\label{zero-00}
Suppose that  $\{x^k\}$, $\{v^k\}$  and $\{w^k\}$ are the sequences generated by Algorithm SEML. If the sequences are bounded, then
$$\lim_{k\rightarrow \infty}\|v^{k+1}-v^k\|=\lim_{k\rightarrow \infty}\|v^k-x^k\|=\lim_{k\rightarrow \infty}\|x^k-w^k\|=0.$$
\end{lemma}

\begin{proof}
The definition of 
$v^{k+1}$ implies that $v^{k+1}\in N_k$. Therefore we have $\langle  v^{k+1}-v^k, g'(v^0)-g'(v^k)\rangle \leq0$ by proposition \ref{lem-proj} which implies that
$$ D_g(v^k,v^0)+ D_g(v^{k+1},v^k)- D_g(v^{k+1},v^0)\leq 0.$$
 Hence $ D_g(v^k,v^0)\leq  D_g(v^{k+1},v^0)$. So, the sequence $\{ D_g(v^k,v^0)\}$ is non-decreasing. Since $\{v^k\}$ is bounded,  $\lim_{k\to\infty} D_g(v^k,v^0)$ exists.
We also have
$$ D_g(v^{k+1},v^k)\leq D_g(v^{k+1},v^0)- D_g(v^k,v^0).$$
Passing to the limit in the above inequality as $k\rightarrow \infty$, we get
$$\lim_{k\rightarrow \infty} D_g(v^{k+1},v^k)=0.$$
Now, by Proposition \ref{zero}, we have
\begin{equation}\label{shom1}
\lim_{k\rightarrow \infty}\|v^{k+1}-v^k\|=0.
\end{equation}
Since $v^{k+1} \in M_k$, from the definition of $M_k$, we have
\begin{equation}\label{pipo}
 \gamma_k D_g(v^k,x^k)\leq \langle v^k- v^{k+1},g'(v^k)-g'(x^k)\rangle.
 \end{equation}
 Therefore, by the Cauchy-Schwarz inequality, we have 
 \begin{equation}\label{}
 \gamma_k D_g(v^k,x^k)\leq \|v^k- v^{k+1}\|\|g'(v^k)-g'(x^k)\|.
 \end{equation} 
 Note that $g'$ is bounded on bounded subsets of $E$ by H3 and Proposition \ref{uni-bound}. Now since $\{x^k\}$ and $\{v^k\}$ are bounded, $\lim_{k\rightarrow \infty}\|v^{k+1}-v^k\|=0$ and $\gamma_k\geq\varepsilon>0$, we get
 $$\lim_{k\rightarrow \infty}D_g(v^k,x^k)=0.$$
 Therefore Proposition \ref{zero} implies that 
 \begin{equation}\label{shom1.5}
  \lim_{k\rightarrow \infty}\|v^k-x^k\|=0.
 \end{equation} 
In the sequel, note that
$$\|v^{k+1}-x^k\|\leq \|v^{k+1}-v^k\|+\|v^k- x^k\|,$$
hence, by \eqref{shom1} and \eqref{shom1.5}, we have  
\begin{equation}\label{shom2}
\lim_{k\rightarrow \infty}\|v^{k+1}-x^k\|=0.
\end{equation}  
On the other hand, since $v^{k+1} \in L_k$, from the definition of $L_k$, we have
\begin{equation}\label{pipo2}
\gamma_k D_g(x^k,w^k)\leq \langle x^k- v^{k+1},g'(x^k)-g'(w^k)\rangle.
 \end{equation} 
Again, by the Cauchy-Schwarz inequality, we have 
\begin{equation}\label{pipo3}
\gamma_k D_g(x^k,w^k)\leq \|x^k- v^{k+1}\|\|g'(x^k)-g'(w^k)\|.
\end{equation} 
Since the sequences $\{x^k\}$ and $\{w^k\}$ are bounded and $\lim_{k\rightarrow \infty}\|v^{k+1}-x^k\|=0$, a similar argument as above shows that  $\lim_{k\rightarrow \infty} D_g(x^k,w^k)=0$. Again  Proposition \ref{zero} implies that $$\lim_{k\rightarrow \infty}\|x^k-w^k\|=0.$$
 \end{proof}

\begin{proposition}\label{pv-2} 
Let $\{x^k\}$, $\{y^k\}$, $\{z^k\}$ and $\{w^k\}$ be the
sequences generated by Algorithm SEML.
If the sequences $\{x^k\}$ and $\{w^k\}$ are bounded and the algorithm does not have finite termination, then\\
i) the sequence $\{z^k\}$ is bounded,\\
ii) there exists a positive sequence $\{\varepsilon_k\}$ such that $\varepsilon_k\to 0$ and $-f(y^k,x^k)+\varepsilon_k e^k \in {\rm int}(C)$ for all $k$.
\end{proposition}

\begin{proof}
i) Since $x^k\in T(v^k)$ by \eqref{ka-1s}, we conclude from \eqref{ka-2s} 
and Theorem \ref{theo-21}, that there exists $c^k\in C^+\setminus\{0\}$ such that
\begin{align} \label{e100}
\nonumber
&\beta_k\langle f(x^k,z^k),c^k\rangle+g(z^k)\langle e^k,c^k\rangle-\langle z^k,g'(x^k)\rangle \langle e^k,c^k\rangle\\
\nonumber
&\le\beta_k\langle f(x^k,x^k),c^k\rangle+g(x^k)\langle e^k,c^k\rangle-\langle x^k,g'(x^k)\rangle\langle e^k,c^k\rangle\\
&=g(x^k)\langle e^k,c^k\rangle-\langle x^k,g'(x^k)\rangle\langle e^k,c^k\rangle,
\end{align}
using B1 in the equality.
From \eqref{e100}, we get 
\begin{align}\label{e101}
\nonumber
g(z^k)\langle e^k,c^k\rangle
&\le -\beta_k\langle f(x^k,z^k),c^k\rangle+\langle z^k,g'(x^k)\rangle \langle e^k,c^k\rangle + g(x^k)\langle e^k,c^k\rangle-\langle x^k,g'(x^k)\rangle\langle e^k,c^k\rangle\\
&= 
-\beta_k\langle f(x^k,z^k),c^k\rangle+\Big(\langle z^k,g'(x^k)\rangle + g(x^k)-\langle x^k,g'(x^k)\rangle\Big)\langle e^k,c^k\rangle.
\end{align}
Take now $u^k_*\in \partial \langle f(x^k, \cdot), c^k\rangle(x^k)$ 
and define  $u^k=\langle e^k,c^k\rangle^{-1} u^k_*$. 
By the definition of $\partial \langle f(x^k, \cdot), c^k\rangle$  evaluated
at $x^k$, we have
\begin{equation}\label{e102}
\langle y-x^k,u^k\rangle \langle e^k,c^k\rangle\le 
\langle f(x^k,y),c^k\rangle-\langle f(x^k,x^k),c^k\rangle=\langle f(x^k,y),c^k\rangle
\end{equation}

Let $B_1(x^k)$ be the closed ball of radius one centered at $x^k$.
Since $f$ is bounded on bounded sets by B2, and $\{x^k\}$ is bounded, there is $M>0$ such that $\lV f(x^k,y)\rV<M$ for all $k$ and for all $y\in B_1(x^k)$. Now without loss of generality, we can assume $\lV c^k\rV=1$ for all $k$. 
Then we have
\begin{equation}\label{e102-a1}
\|u^k\|\langle e^k,c^k\rangle=\sup_{y\in B_1(x^k)} \langle y-x^k,u^k\rangle\langle e^k,c^k\rangle \le \sup_{y\in B_1(x^k)}\langle f(x^k,y),c^k\rangle\leq M.
\end{equation}
Since the sequence $\{e^k\}$ converges strongly to a point in ${\rm int}(C)$ and $\{c^k\}$ is bounded, this implies that $\liminf_{k\to\infty}\langle e^k,c^k\rangle>0$. Therefore \eqref{e102-a1} shows that $\{u^k\}$ is bounded. Now from \eqref{e102}, we have
\begin{equation}\label{e102-b}
\langle z^k-x^k,u^k\rangle \langle e^k,c^k\rangle\le 
\langle f(x^k,z^k),c^k\rangle
\end{equation}
Combining \eqref{e101} and \eqref{e102-b}, we get, after dividing by $\langle e^k,c^k\rangle$, 
\begin{align}\label{e103}
\nonumber
g(z^k)&\le \beta_k\langle x^k-z^k,u^k\rangle+\langle z^k,g'(x^k)\rangle + g(x^k)-\langle x^k,g'(x^k)\rangle\\
&\le \tilde\beta\lV x^k-z^k\rV \lV u^k\rV +\lV z^k\rV \lV g'(x^k)\rV + g(x^k)+\lV x^k\rV \lV g'(x^k)\rV. 
\end{align}
Since $\{x^k\}$ is bounded, and $g'$ is uniformly continuous on bounded subsets of $E$ by H3, we get $g(x^k)$ and $g'(x^k)$ are bounded by Proposition \ref{uni-bound}. Therefore H4 and \eqref{e103} imply that the sequence $\{z^k\}$ is bounded.

ii) Note that $\{x^k\}$ and $\{w^k\}$ are bounded by hypothesis, and $\{y^k\}$ is bounded by item (i) and  \eqref{ka-7s}. 
 Also since $f(\cdot,\cdot)$ is uniformly continuous on bounded sets by B2, and $\lim_{k\rightarrow +\infty} \lV w^k-x^k\rV=0$ by Lemma \ref{zero-00}, we conclude that
\begin{equation}\label{e105}
\lim_{k\rightarrow +\infty}\lV f(y^k,x^k)-
f(y^k,w^k)\rV =0.
\end{equation}
Also note that $C$ is a closed and convex cone, and $f$ is bounded on bounded sets because $f$ is uniformly continuous on bounded sets. Thus $\{ f(y^k,x^k)\}$ and $\{ f(y^k,w^k)\}$ are bounded. Now since 
 $w^k\in H_k$ by \eqref{ka-8s},  $f(y^k,w^k)\in -C$ for all $k$ by \eqref{ka-8sh}.  Now this fact, together with \eqref{e105}, easily imply that 
there exists a positive sequence $\{\varepsilon_k\}$ such that $\varepsilon_k\to 0$ and $-f(y^k,x^k)+\varepsilon_k e^k \in {\rm int}(C)$ for all $k$.
\end{proof}

\begin{proposition}\label{fix-Ks}
 Assume that $f$ is a vector valued bifunction, $T(\cdot)$ is a multivalued  mapping  and the assumptions B1--B4 are satisfied.\\
i) If there exists a subsequence $\{x^{k_n}\}$ of $\{x^k\}$ such that $x^{k_n}\rightharpoonup p$,  then $p\in K_{\infty}\cap{\rm Fix}(T)$,  where $K_{\infty}=\cap_{k=0}^{\infty}K_k$.\\
ii) $K_{\infty}\cap{\rm Fix}(T)\subset L_k \cap M_k\cap N_k$ for all $k$.
\end{proposition}
\begin{proof} 
i) We first prove that
$p\in {\rm Fix}(T)$. Note that we have $\lim_{n\rightarrow \infty} \lV v^{k_n}-x^{k_n}\rV=0$ by Lemma \ref{zero-00}, where for each $n$, $x^{k_n}$ is the Bregman projection of $v^{k_n}$ onto $T(v^{k_n})$. Therefore we have 
$\lim_{n\rightarrow \infty}d(v^{k_n},T(v^{k_n}))=0$. Now since $T$ is demiclosed, $p\in T(p)$, i.e. 
$p$ is a fixed point of $T(\cdot)$. Now we prove that $p\in K_{\infty}$. Since $K_{\infty}=\cap_{k=0}^{\infty}K_k$, it suffices to prove that $p\in K_k$ for all
 $k$. Note that the sequence $\{K_k\}$ is nonincreasing, now let $m$ be a fixed integer, hence there is $j>m$ such that for all $n\geq j$ we have
$$w^{k_n}\in K_{k_n}\subset K_m, \ \ \ \ \ \ \ \forall n\geq j,$$
where $w^{k_n}=\Pi^g_{K_{k_n}}(x^{k_n})$. Now, since $\lim_{n\rightarrow \infty} \lV w^{k_n}-x^{k_n}\rV=0$ by Lemma \ref{zero-00}, we have $w^{k_n}\rightharpoonup p$. Consequently, since $K_m$ is closed and convex, we conclude that $p\in K_m$ for all $m$, and hence $$p\in \cap_{k=0}^{\infty}K_k=K_{\infty}.$$
ii) The proof is similar to the proof of Lemma \ref{xn}. It suffices to replace $DS(f,T)$ by $K_{\infty}\cap{\rm Fix}(T)$.
\end{proof}
\begin{remark}\label{remak-f}
It is easy to see that $DS(f,T)\subset K_{\infty}\cap{\rm Fix}(T)$. Also, since $T$ is quasi $ D_g$-nonexpansive,  $K_{\infty}\cap{\rm Fix}(T)$ is closed and convex.
\end{remark}

 In the following proposition, we prove that the sequences $\{v^k\}$ and $\{x^k\}$ generated by Algorithm SEML converge strongly  to
an element of $K_{\infty}\cap{\rm Fix}(T)$.
\begin{proposition} \label{main-seq}
Assume that $f$ is a vector valued bifunction, $T(\cdot)$ is a multivalued  mapping  and the assumptions B1--B4 are satisfied. If the sequences $\{x^k\}$ and $\{v^k\}$ generated by Algorithm SEML are bounded, then the sequences $\{x^k\}$ and $\{v^k\}$ are strongly convergent to
$\bar{x}= \Pi^g_{K_{\infty}\cap{\rm Fix}(T)} (v^0)$.
\end{proposition}
\begin{proof}
Assume that $p$ is any weak limit point of the sequence $\{x^k\}$. Then, there exists a
subsequence $\{x^{k_n}\}$ of $\{x^k\}$ such that $x^{k_n} \rightharpoonup p$ as $n\rightarrow \infty$.
 Note that  Proposition \ref{fix-Ks} shows that $p\in  K_{\infty}\cap{\rm Fix}(T)$ and hence $ K_{\infty}\cap{\rm Fix}(T)\not=\emptyset$. Also, $K_{\infty}\cap{\rm Fix}(K)$ is closed and convex by Remark \ref{remak-f}, therefore  $\bar{x}= \Pi^g_{K_{\infty}\cap{\rm Fix}(T)} (v^0)$ is well defined. In the sequel, we first prove the weak convergence of the sequence $\{x^k\}$.
Then we show that $x^k \rightarrow \bar{x}= \Pi^g_{K_{\infty}\cap{\rm Fix}(T)} (v^0)$. 
From the definition of $N_k$, we have $v^k=\Pi^g_{N_k}(v^0)$. Since $K_{\infty}\cap{\rm Fix}(T)\subset N_k$ by Proposition \ref{fix-Ks} (ii), and  $\Pi^g_{N_k}$ is the Bregman projection map onto $N_k$,  for  $\bar{x}\in K_{\infty}\cap{\rm Fix}(T) \subset N_k$, we have $\langle \bar{x}-v^k, g'(v^0)-g'(v^k) \rangle \leq0$ by Proposition \ref{lem-proj}. This implies that $ D_g(v^k,v^0)\leq  D_g(\bar{x},v^0)$.
Therefore we have
\begin{equation}\label{hhh}
g(v^k)-g(v^0)-\langle v^k-v^0, g'(v^0) \rangle \leq  D_g(\bar{x},v^0).
\end{equation}
Since $v^{k_n} \rightharpoonup p$ by Lemma \ref{zero-00},  by the weak lower semicontinuity of $g$ and replacing $k$ by $k_n$ in \eqref{hhh}, letting $n\to\infty$, we get
$$ D_g(p,v^0)=g(p)-g(v^0)-\langle p-v^0, g'(v^0) \rangle\leq \liminf_{n\to\infty} \Big(g(v^{k_n})-g(v^0)-\langle v^{k_n}-v^0, g'(v^0) \rangle\Big)\leq  D_g(\bar{x},v^0).$$
From the definition of $\bar{x}$ and since $p\in  K_{\infty}\cap{\rm Fix}(T)$, we get $\bar{x}=p$, i.e. $x^{k_n} \rightharpoonup \bar{x}$. Hence every weakly convergent subsequence of $\{x^k\}$ converges weakly to $\bar{x}$.
This shows that $x^k\rightharpoonup \bar{x}$, and therefore $v^k \rightharpoonup \bar{x}$. Taking liminf and limsup in \eqref{hhh}, we get $\lim_{k\rightarrow \infty} g(v^k)=g(\bar{x})$. This implies that
$$\lim_{n\rightarrow \infty}  D_g(v^k,\bar{x})=\lim_{n\rightarrow \infty}\Big( g(v^k)-g(\bar{x})-\langle v^k-\bar{x}, g'(\bar{x}) \rangle\Big)=0.$$
Therefore by proposition \ref{zero}, we have $v^k \rightarrow \bar{x}= \Pi^g_{K_{\infty}\cap{\rm Fix}(T)} (v^0)$. Now since, $\lim_{k\rightarrow \infty}\|v^k-x^k\|=0$ by Lemma \ref{zero-00}, we get $x^k \rightarrow \bar{x}= \Pi^g_{K_{\infty}\cap{\rm Fix}(T)}(v^0)$. 
\end{proof}

\begin{proposition}\label{asy-sol-2} 
Let $\{x^k\}$ and $\{z^k\}$ be the
sequences generated by Algorithm SEML.
If $\{x^{k_i}\}$ is a subsequence of $\{x^k\}$ satisfying $\lim_{i\rightarrow+\infty} D_g(z^{k_i},x^{k_i})=0$, 
 then $\bar{x}\in S(f,T)$ where $\bar{x}= \Pi^g_{K_{\infty}\cap{\rm Fix}(T)} (v^0)$ is the strong limit of $\{x^k\}$.
\end{proposition}
\begin{proof}
Since $\lim_{i\rightarrow+\infty} D_g(z^{k_i},x^{k_i})=0$, Proposition \ref{zero} implies that 
\begin{equation}\label{e107}
\lim_{i\rightarrow+\infty}\lV z^{k_i}-x^{k_i}\rV =0.
\end{equation}
Now since $g'$ is uniformly continuous on bounded sebsets of $E$ by H3, we get from \eqref{e107}, 
\begin{equation}\label{e111}
\lim_{i\rightarrow+\infty}\lV g'(z^{k_i})-g'(x^{k_i})\rV =0. 
\end{equation}
On the other hand, since $\{x^k\}$ and $\{z^k\}$ are bounded, 
   from B2 and \eqref{e107}, we obtain
\begin{equation}\label{e109}
\lim_{i\rightarrow +\infty} f(x^{k_i},z^{k_i})=0.
\end{equation}
Note that $x^{k_i}\rightarrow \bar{x}$ and $v^{k_i}\rightarrow \bar{x}$ by Proposition \ref{main-seq}. Now take any $y\in T(\bar{x})$, since $T$ is  lower semicontinuous at $\bar{x}\in K$, there is a sequence $\{\xi^{k_i}\}$ such that $\xi^{k_i}\in T(v^{k_i})$ and
$\xi^{k_i}\rightarrow y$. By Corollary \ref{c1}, we have  
$$ 
\langle \xi^{k_i}-z^{k_i},g'(x^{k_i})-g'(z^{k_i})\rangle\langle e^{k_i},c^{k_i}\rangle\leq
\beta_{k_i}\left[\langle f(x^{k_i},\xi^{k_i}),c^{k_i}\rangle-\langle f(x^{k_i},z^{k_i}),c^{k_i}\rangle\right],
$$
which implies that
\begin{equation}\label{vay-1}
-\beta_{k_i}^{-1}\lV \xi^{k_i}-z^{k_i}\rV \lV g'(x^{k_i})-g'(z^{k_i})\rV\langle e^{k_i},c^{k_i}\rangle\le
\langle f(x^{k_i},\xi^{k_i}),c^{k_i}\rangle-\langle f(x^{k_i},z^{k_i}),c^{k_i}\rangle.
\end{equation}
Without loss of generality, we can assume that $c^{k_i}\rightharpoonup c^*\in C^+\setminus\{0\}$.
Taking the limit from \eqref{vay-1}, and using \eqref{e111} and \eqref{e109}, we conclude that 

\begin{equation}\label{vay-2}
0\le \lim_{i\to\infty}\langle f(x^{k_i},\xi^{k_i}),c^{k_i}\rangle=\langle f(\bar{x},y),c^*\rangle. 
\end{equation} 
Note that $c^*\in C^+\setminus\{0\}$, therefore we have 
$$
f(\bar{x},y)\not\in -{\rm int}(C).
$$
 Since $y\in K(\bar{x})$ is arbitrary,  $\bar{x}\in S(f,T)$.
\end{proof}

\begin{proposition}\label{alpha-0-2}
If a subsequence $\{\alpha_{k_i}\}$ of $\{\alpha_k\}$ as defined in \eqref{ka-6s} converges to $0$,  then $\bar{x}\in S(f,T)$ where $\bar{x}= \Pi^g_{K_{\infty}\cap{\rm Fix}(T)} (v^0)$ is the strong limit of $\{x^k\}$.
\end{proposition}
\begin{proof}
For proving the result, we will use Proposition \ref{asy-sol-2}. Thus, we must show that 
$$
\lim_{i\rightarrow+\infty} D_g(z^{k_i},x^{k_i})=0.
$$
For the sake of contradiction, and without loss of generality, let us assume that
\begin{equation}\label{e112}
{\rm liminf}_{i\rightarrow+\infty} D_g(z^{k_i},x^{k_i})\ge\eta>0, 
\end{equation}
taking into account the nonnegativity of $D_g(\cdot,\cdot)$.
Define
\begin{equation}\label{e201}
\hat{y}^i=\frac{\alpha_{k_i}}{\theta}z^{k_i} +(1-\frac{\alpha_{k_i}}{\theta})x^{k_i},
\end{equation}
where $\alpha_{k_i}=\theta^{\ell(k_i)}$ by \eqref{ka-6s}.
Therefore we have
\begin{equation}\label{e202}
\hat{y}^i-x^{k_i}=\frac{\alpha_{k_i}}{\theta}(z^{k_i}-x^{k_i}).
\end{equation}
Note that  $\lim_{i\rightarrow+\infty}\alpha_{k_i}=0$, hence $\ell(k_i)> 1$ for large enough $i$.
Also, in view of \eqref{e201}, we have that 
$\hat y^i=y^{\ell(k_i)-1}$ in the inner loop of the linesearch for determining $\alpha_{k_i}$, 
i.e., in \eqref{ka-5s}. Since $\ell(k_i)$ is the first integer for which the exclusion in \eqref{ka-4s} holds,
such exclusion does not hold for $\ell(k_i)-1$. i.e., we have 
\begin{equation}\label{e113}
 -\beta_{k_i}f(\hat{y}^i,x^{k_i})+\beta_{k_i}f(\hat{y}^i,z^{k_i})+ \delta D_g(z^{k_i},x^{k_i})e^{k_i}\in {\rm int}(C)
\end{equation}
for large enough $i$. On the other hand, 
since $\lim_{i\rightarrow+\infty}\alpha_{k_i}=0$ by hypothesis, and $\{z^{k_i}-x^{k_i}\}$ is bounded 
by  Lemma \ref{bound-vk} and Proposition \ref{pv-2}(i),  it follows from \eqref{e202} that 
\begin{equation}\label{we-4}
\lim_{i\rightarrow+\infty}\lV\hat{y}^i-x^{k_i}\rV=0.
\end{equation}
Since $f(\cdot,\cdot)$ is uniformly continuous on bounded sets by B2, \eqref{e113} and \eqref{we-4} imply that
\begin{equation}\label{we-4.1}
 -\beta_{k_i}f(x^{k_i},x^{k_i})+\beta_{k_i}f(x^{k_i},z^{k_i})+ \delta D_g(z^{k_i},x^{k_i})e^{k_i}\in C
\end{equation}
for large enough $i$.
Since $\delta$ belongs to $(0,1)$, it follows from \eqref{we-4.1} that
\begin{equation}\label{we-4.2}
 \beta_{k_i}f(x^{k_i},z^{k_i})+  D_g(z^{k_i},x^{k_i})e^{k_i}\in {\rm int}(C).
\end{equation}
Take now $y=x^{k_i}$ in Corollary \ref{c1}, then we have  
\begin{align}
\nonumber
\left[ D_g(z^{k_i},x^{k_i})+ D_g(x^{k_i},z^{k_i})\right]\langle e^{k_i},c^{k_i}\rangle
&=\langle x^{k_i}-z^{k_i},g'(x^{k_i})-g'(z^{k_i})\rangle\langle e^{k_i},c^{k_i}\rangle\\
&\leq
\beta_{k_i}\left[\langle f(x^{k_i},x^{k_i}),c^{k_i}\rangle-\langle f(x^{k_i},z^{k_i}),c^{k_i}\rangle\right],
\end{align}
which implies that
\begin{equation}\label{we-4.3}
\left[ D_g(z^{k_i},x^{k_i})+ D_g(x^{k_i},z^{k_i})\right]\langle e^{k_i},c^{k_i}\rangle+
\beta_{k_i}\langle f(x^{k_i},z^{k_i}),c^{k_i}\rangle\leq 0.
\end{equation}
Since $c^{k_i}\in C^+\setminus\{0\}$, we have 
\begin{equation}\label{we-4.4}
D_g(z^{k_i},x^{k_i})e^{k_i}+ D_g(x^{k_i},z^{k_i})e^{k_i}+ \beta_{k_i}f(x^{k_i},z^{k_i})\not\in {\rm int}(C).
\end{equation}
Note that $ D_g(x^{k_i},z^{k_i})>0$, hence \eqref{we-4.4}  contradicts \eqref{we-4.2}, thus establishing the result.

\end{proof}

We now complete the paper by giving the proof of Theorem \ref{main-strong}. 

\begin{flushleft}
\textbf{Proof of Theorem \ref{main-strong}:}
\end{flushleft}

\begin{proof} 
Note that $x^k \rightarrow \bar{x}= \Pi^g_{K_{\infty}\cap{\rm Fix}(T)}(v^0)$ by Proposition \ref{main-seq}.  We consider two cases related to the behavior  of $\{\alpha_k\}$.
First assume that there exists a subsequence $\{\alpha_{k_i}\}$ of $\{\alpha_k\}$ which converges to $0$. In this case, the result is obtained by Proposition \ref{alpha-0-2}, i.e. we get that $\bar{x}\in S(f,T)$.

Now we take a subsequence $\{\alpha_{k_i}\}$ of $\{\alpha_k\}$ bounded away from zero, say greater or equal to $\eta$
for large enough $i$. It follows
from \eqref{ka-4s} and \eqref{ka-7s} that
\begin{equation}\label{we-5}
-\beta_{k_i}f(y^{k_i},x^{k_i})+\beta_{k_i}f(y^{k_i},z^{k_i})+ \delta D_g(z^{k_i},x^{k_i})e^{k_i}\not\in {\rm int}(C).
\end{equation}
Note that, since $\alpha_{k_i}\le 1$ by \eqref{ka-6s}, we get, in view of B1 and B3, 
\begin{equation}\label{e117}
0= f(y^{k_i},y^{k_i})\preceq \alpha_{k_i} f(y^{k_i},z^{k_i})+(1-\alpha_{k_i}) f(y^{k_i},x^{k_i})\in C.
\end{equation}
Hence we have
\begin{equation}\label{e117.5}
-\beta_{k_i} f(y^{k_i},z^{k_i})+\frac{-\beta_{k_i}(1-\alpha_{k_i})}{\alpha_{k_i}} f(y^{k_i},x^{k_i})\in -C.
\end{equation}
Summing up \eqref{we-5} and \eqref{e117.5}, we have 
\begin{equation}\label{we-5-2}
\frac{-\beta_{k_i}}{\alpha_{k_i}}f(y^{k_i},x^{k_i})+ \delta D_g(z^{k_i},x^{k_i})e^{k_i}\not\in {\rm int}(C).
\end{equation}
We claim that 
$
\lim_{i\rightarrow+\infty} D_g(z^{k_i},x^{k_i})=0
$. If our claim were false, then without loss of generality, there would exist $\lambda>0$ such that
$$
\lim_{i\rightarrow+\infty} D_g(z^{k_i},x^{k_i})>\lambda>0.
$$ 
Therefore \eqref{we-5-2} implies that 
\begin{equation}\label{we-5-201}
-f(y^{k_i},x^{k_i})+ \frac{\alpha_{k_i}\delta\lambda}{\beta_{k_i}}e^{k_i}\not\in {\rm int}(C)
\end{equation}
for large enough $i$. Since the positive sequence $\{\frac{\alpha_{k_i}\delta\lambda}{\beta_{k_i}}\}$ is  bounded away from zero,  \eqref{we-5-201} contradicts Proposition \ref{pv-2} (ii). Therefore
$$
\lim_{i\rightarrow+\infty} D_g(z^{k_i},x^{k_i})=0.
$$  
Now, we invoke Proposition \ref{asy-sol-2} in order to get  $\bar{x}\in S(f,T)$. We have shown that the limit $\bar{x}$ of $\{x^k\}$
belongs to $S(f,T)$ both when the corresponding stepsizes $\{\alpha_k\}$ 
either approach zero or remain bounded away from zero, establishing the claim.
\end{proof}

\section{Examples and numerical experiments }\label{s4}

In this section, we first give some examples of vector quasi-equilibrium problems in several Banach spaces to which  our main theorem can be applied for finding a solution. Then we present some numerical experiments.

\begin{example} 
Suppose that $E=\ell^p=\Big\{\xi=(\xi_1, \xi_2 ,\xi_3, \cdots): \|\xi\|_p=(\sum_{i=1}^{\infty}|\xi_i|^p)^{\frac{1}{p}}< \infty\Big\}$ for $1<p<\infty$,
  $K=\Big\{\xi=(\xi_1, \xi_2 ,\xi_3, \cdots)\in \ell^p: \xi_i\geq0, \
  i=1,2,3 \Big\}$  and 
 $C=\Big\{z\in \mathbb{R}^3: z_i\geq0, i=1,2,3\Big\}$.  
 We define the vector valued bifunction $f:E\times E \rightarrow\mathbb{R}^3$ as $$f(x,y)=\langle  y-x, J(x-A(x))\rangle (3x_1x_2+1, x_3^2+x_2+2, 7x_2^2+4x_1x_3+1)$$ where the map $A: E\to E$ is defined by
$$A(x)=(x_1^2+x_1-9, 3x_2-5,  x_3^3+x_3-8 , -x_4, -x_5, \cdots ).$$  We also define $T(\cdot):K\to {\cal P}(K)$ by
$T(x)=\Big\{\xi\in K: \|\xi\|_p\leq \|x\|_p \Big\}$. If $x$ is an equilibrium point of VQEP$(f, T)$, then we have
\begin{equation}
\label{VEP-s4e1}
f(x,y)=\langle  y-x, J(x-A(x))\rangle (3x_1x_2+1, x_3^2+x_2+2, 7x_2^2+4x_1x_3+1)\not\in -{\rm int}(C)
\end{equation}
for all $y\in T(x)$. If $A(x)\in T(x)$, then we get
\begin{equation}
\label{VEP-s4e2}
\langle  A(x)-x, J(x-A(x))\rangle (3x_1x_2+1, x_3^2+x_2+2, 7x_2^2+4x_1x_3+1)\not\in -{\rm int}(C).
\end{equation}
 Therefore there 
exists $c^*\in C^+\setminus\{0\}$ such that
\begin{equation}
\label{VEP-s4e3}
\langle  A(x)-x, J(x-A(x))\rangle \langle (3x_1x_2+1, x_3^2+x_2+2, 7x_2^2+4x_1x_3+1), c^*\rangle \geq0.
\end{equation}
Since $x\in K$, we have
 $(3x_1x_2+1, x_3^2+x_2+2, 7x_2^2+4x_1x_3+1)\in {\rm int}(C)$. Thus we get
  $$\langle (3x_1x_2+1, x_3^2+x_2+2, 7x_2^2+4x_1x_3+1), c^*\rangle >0.$$
  Therefore \eqref{VEP-s4e3} implies that
 \begin{equation}
\label{VEP-s4e4}
-\| A(x)-x \|_p^2 =\langle  A(x)-x, J(x-A(x))\rangle\geq 0.
\end{equation}
Now \eqref{VEP-s4e4} shows that
$$-(|x_1^2-9|^p+|2x_2-5|^p+| x_3^3-8|^p+ \sum_{i=4}^{\infty}|2x_i|^p)^\frac{2}{p} \geq 0 .$$ Hence, we get $x_1=\pm 3$, $x_2=\frac{5}{2}$, $x_3=2$ and 
$x_i=0$ for all $i\geq 4$. On the other hand, since $x\in K$,  we conclude that $x_1=3$. Therefore $x=(3, \frac{5}{2}, 2,  0, 0, \cdots)$ is a solution of VQEP$(f,T)$. Note that $(0, 0, 0, \cdots)$ is also another solution of the problem.

Moreover, it is obvious that the assumptions B1-B3 are satisfied and $T$ is a multivalued mapping with nonempty, closed and convex values.  Since  $x\in T(x)$ for all $x\in K$, it is easy to see that
$T$ is a demiclosed and  quasi $D_g$-nonexpansive mapping with $g(\cdot)=\|\cdot\|_p^2$. Now we show that $T$ is lower semicontinuous at each $\bar{p}\in K$.   Suppose that $p^k\rightarrow \bar{p}$ and $\bar{q}\in T(\bar{p})$. Define $q^k:=\Pi^g_{T(p^k)}(\bar{q})$  where $\Pi^g$ is the Bregman projection. Hence we have  $q^k\in T(p^k)$ and
$q^k\to \bar{q}$. This means $T$ is lower semicontinuous at each $\bar{p}\in K$. Therefore the assumption B4 is satisfied. Now, if $DS(f,T)\neq\emptyset$ or  the sequences generated by Algorithm SEML are bounded,
then Theorem \ref{main-strong} ensures that  the sequence $\{v^k\}$
converges strongly to a solution of the problem.

\end{example} 
We continue with another example of vector quasi-equilibrium problems in $L^p$ spaces with $p=2$ (the set of square integrable functions) to which our main result can be applied.
\begin{example} 
Suppose that $E= L^2([a,b])$ with the inner product
$$\langle x, y \rangle:=\int_a^b x(t)y(t) dt, \ \ \ \ \ \ \ \ \forall \ x, y \in E $$
and the induced norm
$$\|x\|:=\Big (\int_a^b |x(t)|^2 dt\Big )^\frac{1}{2}, \ \ \ \ \ \ \ \ \forall \ x \in E. $$
Let $I$ be a totally ordered finite set, suppose that $Y=\ell^2(I)$, and
let $\{\varepsilon_i\}_{i\in I}$ be the standard unit vectors. Then
$C=\ell_+^2(I)=\Big\{\xi=(\xi_i)_{i\in I}: \xi_i\geq 0, \ \forall i\in I\Big\}$ is a closed convex pointed cone with nonempty interior.
We define the vector valued bifunction $f:E\times E \rightarrow Y$ as
$$f(x,y)=\sum_{i\in I}\phi_i(x)(\psi_i(y)-\psi_i(x))\varepsilon_i$$   
where $\phi_i: E \to \mathbb{R}^+$  is uniformly continuous on bounded sets for all $i\in I$, and also $\psi_i:E\to \mathbb{R}$
 is convex and uniformly continuous on bounded sets for all $i\in I$.   
Take $K=E$ and define $T(\cdot):K\to {\cal P}(K)$ by  $T(x)=\Big\{\frac{t}{3}x+(1-t)z: 0\leq t \leq 1, z\in B(x,\frac{1}{3}\|x\|) \Big\}$ for each $x\in K$, where $B(x,\frac{1}{3}\|x\|)$ denotes the closed ball of radius $\frac{1}{3}\|x\|$ centered at $x$.  
It is easy to see that  $f$ satisfies B1--B3. Now we show that B4 is satisfied. Note that $T$ is a multivalued mapping with nonempty, closed and convex values. It is obvious that $T$ is demiclosed, because  $x\in T(x)$ for all $x\in K$.  Take $p\in {\rm Fix}(T)$ and $x\in K$, then we have $\|p-P_{T(x)}(x)\|= \|p-x\|$ where $P$ is the metric projection. This shows that $T$ is quasi $D_g$-nonexpansive with $g(\cdot)=\|\cdot\|^2$. Now, suppose that $x^k\rightarrow \bar{x}$ and $\bar{y}\in T(\bar{x})$. Therefore there exist $t\in [0,1]$ and $\bar{z}\in B(\bar{x},\frac{1}{3}\|\bar{x}\|)$ such that $\bar{y}=\frac{t}{3}\bar{x}+(1-t)\bar{z}$. Define $z^k:= P_{B(x^k,\frac{1}{3}\|x^k\|)}(\bar{z})$ where $P$ is the metric projection. Thus we get $z^k\to \bar{z}$. Now we define
 $y^k=\frac{t}{3}x^k+(1-t)z^k$.
 Then we have $y^k\in T(x^k)$ and
$y^k\to \bar{y}$. This shows that $T$ is lower semicontinuous at each $\bar{x}\in K$. Now if $\bigcap_{i\in I}{\rm argmin}\psi_i\not=\emptyset$ or 
 the sequences generated by Algorithm SEML are bounded, then 
 Theorem \ref{main-strong} ensures that the sequence $\{v^k\}$ converges strongly to a solution of VQEP$(f,T)$.

\end{example}

\begin{flushleft}
\textbf{The (Generalized)  Nash Equilibrium Problem:}
\end{flushleft}
Suppose that $I = \{1, 2, \cdots ,n\}$ is a finite index set which denotes the set of players. Let $E_i$ be a Banach space where $i\in I$, and the strategy set $K_i$ is subset of $E_i$ for the $i$-th player. Note that $E := E_1 \times E_2 \times \cdots \times E_n$ is a Banach space, and the set $K := K_1 \times K_2 \times \cdots \times K_n$ is a subset of the
Banach space  $E$. Let $\varphi_i : K \to \mathbb{R}$ be a payoff function which shows the loss of each player where $i\in I$. Also, $\varphi_i$   depends on the strategies of all the player for any $i\in I$. The Nash equilibrium problem corresponding
to $\{\varphi_i\}_{i\in I}$ and $\{K_i\}_{i\in I}$ is to find $x=(x_1, x_2, \cdots, x_n)\in K$ such that 
$$\varphi_i(x)\leq \varphi_i(x_1, \cdots , x_{i-1}, y_i , x_{i+1}, \cdots , x_n),$$ 
for all $i\in I$ and all $y_i\in K_i$ .
The point $x$ is a solution of the problem and is called a Nash
equilibrium point. The above inequality implies that each
Nash equilibrium point corresponds to an optimal amount  for minimizing the loss. 
Now, we define $f:K\times K \to \mathbb{R}$ as
$$f(x,y)=\sum_{i=1}^{n}(\varphi_i(x_1, \cdots , x_{i-1}, y_i , x_{i+1}, \cdots , x_n)-\varphi_i(x_1, \cdots , x_{i-1}, x_i , x_{i+1}, \cdots , x_n)),$$
where $x=(x_1, x_2, \cdots, x_n)$ and $y=(y_1, y_2, \cdots, y_n)$.
 So,  $f$ is a bifunction and its corresponding equilibrium problem is to find $x\in K$ such that 
   $$f(x,y)\geq 0, \ \ \ \ \ \ \ \ {\rm for \ all}\ y\in K.$$
 It is easy to see that $x$ is a Nash equilibrium point if and only if $x$ is an equilibrium point of $f$. 
 
Now, we extend this problem to construct a vector quasi-equilibrium problem. Therefore, for any $i\in I$, we extend the payoff function $\varphi_i : K \to \mathbb{R}$ to a finite family of functions $\varphi_{ij} : K \to \mathbb{R}$ showing the loss of the $i$-th player in $m$ areas separately (for example, losses in the areas of finance, energy, time, human resources and etc) where $1\leq j\leq m$. Consider
 $$f_j(x,y)=\sum_{i=1}^{n}(\varphi_{ij}(x_1, \cdots , x_{i-1}, y_i , x_{i+1}, \cdots , x_n)-\varphi_{ij}(x_1, \cdots , x_{i-1}, x_i , x_{i+1}, \cdots , x_n)),$$
 for all $1\leq j\leq m$ and $x=(x_1, x_2, \cdots, x_n)$ and $y=(y_1, y_2, \cdots, y_n)$. We define $f:E\times E\to \re^m$ as
 \begin{equation}\label{vbif-nash}
 f(x,y)=\Big(f_1(x,y), f_2(x,y), \cdots, f_m(x,y)\Big).
 \end{equation}

Consider $C=\Big\{z\in \mathbb{R}^m: z_i\geq0, i=1,2, \cdots, m\Big\}$ which is a closed, convex and pointed cone with nonempty interior.  We also define $T(\cdot):K\to {\cal P}(K)$ by 
\begin{equation}\label{vmap-nash}
T(x)=\prod_{i=1}^mT_i(x)  \ \ \ \ \ \ {\rm and} \ \ \ \ \ \
T_i(x)=\Big\{ y_i\in K_i: h_i(x,y_i)\leq 0\Big\}
\end{equation}
where $h_i: K\times K_i\to \mathbb{R}$ and $h_i(x,x_i)\leq0$ for all $x\in K$ and $1\leq i\leq m$.
Finally, suppose that the optimal amount $x$ to minimize the loss  must belong to $T(x)$.
Therefore,  our problem has been formulated as a vector quasi-equilibrium problem in the Banach space $E$, and the solution of the problem is the vector which minimizes the losses of the $nm$ payoff functions corresponding to the $m$ areas of the problem.

 \begin{example}\label{ex-gnep}
 Consider the vector quasi-equilibrium problem VQEP$(f, T)$  as defined  above and assume that the cost function $\varphi_{ij}$ is convex and uniformly continuous on bounded sets for all $i,j$. Also, for each $i$, suppose that the function  $h_i$ is continuous with respect to both variables and quasi-convex with respect to the second variable.
  It is obvious that the assumptions B1-B3 are satisfied.
Now we show that B4 is satisfied. Note that $T$ is a multivalued mapping with nonempty, closed and convex values because $h_i$ is quasi-convex and continuous with respect to the second variable and $h_i(x,x_i)\leq0$ for all $x\in K$ and $1\leq i\leq m$. Since  $x\in T(x)$ for all $x\in K$, this implies that
$T$ is demiclosed by Definition \ref{defdemic}. On the other hand, by the definition of  quasi $D_g$-nonexpansive mappings, $T$ is quasi $D_g$-nonexpansive because $x=\Pi^g_{T(x)}(x)$  where $\Pi^g$ is the Bregman projection. Now we show that $T$ is lower semicontinuous at each $\bar{x}\in K$.   Suppose that $x^k\rightarrow \bar{x}$ and $\bar{y}\in T(\bar{x})$. Define $y^k:=\Pi^g_{T(x^k)}(\bar{y})$  where $\Pi^g$ is the Bregman projection. It is easy to see that  $y^k\in T(x^k)$ and
$y^k\to \bar{y}$. Therefore $T$ is lower semicontinuous at each $\bar{x}\in K$. Hence the assumption B4 is satisfied. Now, if $DS(f,T)\neq\emptyset$ or  the sequences generated by Algorithm SEML are bounded,
then Theorem \ref{main-strong} ensures that  the sequence $\{v^k\}$
converges strongly to a solution of VQEP$(f,T)$.

\end{example}

 \begin{example} \label{ex-matrix}
 Define the vector valued bifunction
$f:\mathbb{R}^n\times \mathbb{R}^n \rightarrow\mathbb{R}^m$ as
\begin{equation}\label{bimat}f(x,y)=\sum_{i=1}^m\langle A_ix+B_iy+c_i, y-x\rangle \varepsilon_i
\end{equation}
for every $x,y \in\mathbb{R}^n$, where $\{\varepsilon_i\}$ is the standard unit vectors for $\mathbb{R}^m$, and the vector $c_i\in \mathbb{R}^n$, and the matrices $A_i$ and $B_i$ are  square matrices of order $n$ such
that $B_i$ is positive semidefinite for all $1\leq i\leq m$. 
Consider $K= \prod_{i=1}^n[-a_i, a_i]$ where $a_i\in  \mathbb{R}^+ $ and $C=\Big\{z\in \mathbb{R}^m: z_i\geq0, i=1,2, \cdots, m \Big\}$. We also define $T :K\to {\cal P}(K)$ by 
\begin{equation}\label{}
T(x)=\Big \{ z\in K:  \max\{\sum_{i=1}^n x_i, d_1\}\leq \sum_{i=1}^n z_i \  \ {\rm and} \ \ \|z\|\leq \max \{\|x\|, \ d_2\}\Big\}
 \end{equation}
where $ d_1\leq d_2\leq \sum_{i=1}^n a_i$. It is obvious that $f(x,x)=0$ for all $x\in \mathbb{R}^n$, and
 $f$ is $C$-convex with respect to the second variable because $B_i$ is  positive semidefinite for all $1\leq i\leq m$, and 
  $f(\cdot,\cdot)$ is uniformly continuous on bounded sets. Hence  $f$ satisfies B1--B3. Also,  $T$ is a multivalued mapping with nonempty, closed and convex values. If $\{p^k\}$ is a sequence such that  $p^k\to \bar{p}$ and $\lim_{k\to\infty}d(p^k,T(p^k))=0$, it is easy to see that $\bar{p}\in{\rm Fix}(T)$. This implies that $T$ is demiclosed. Take $p\in {\rm Fix}(T)$ and $x\in K$. Then we have $\|p-P_{T(x)}(x)\|\leq \|p-x\|$ where $P$ is the metric projection. This shows that $T$ is quasi $D_g$-nonexpansive with $g(\cdot)=\|\cdot\|^2$.  Now we show that $T$ is lower semicontinuous at each $\bar{x}\in K$.   Suppose that $x^k\rightarrow \bar{x}$ and $\bar{y}\in T(\bar{x})$. Define $y^k:=P_{T(x^k)}(\bar{y})$  where $P$ is the metric projection. It is easy to see that  $y^k\in T(x^k)$ and
$y^k\to \bar{y}$. Therefore $T$ is lower semicontinuous at each $\bar{x}\in K$. Hence the assumption B4 is satisfied. Now since $K$ is bounded,  the sequences generated by Algorithm SEML are bounded. Therefore Theorem \ref{main-strong} ensures that  the sequence $\{v^k\}$
converges strongly to a solution of VQEP$(f,T)$.\\
Now, in order to implement Algorithm SEML in Section \ref{s3} for this example, we  take $n=3$, $m=2$, $d_1=1$, $d_2=3$, $\delta=\frac{1}{1000}$, $\theta=\frac{1}{2}$,  $\beta_k\equiv 1$, $\gamma_k\equiv 1$ and $e^k\equiv (1, 1)$.  We also consider
$$
A_1=\left[ {\begin{array}{ccc}
-1 & 3 & 0  \\
-3 & -2 & 0 \\ 
0 & 0 & -3  \\
\end{array} } \right]
\ \ \  
B_1=\left[ {\begin{array}{ccc}
1 & 0 & -2  \\
0 & 2 & 0 \\ 
2 & 0 & 3  \\
\end{array} } \right]
\ \ \ 
A_2=\left[ {\begin{array}{ccc}
-5 & -1 & 2  \\
1 & -3 & 0 \\ 
-2 & 0 & -2 \\
\end{array} } \right]
\ \ \ 
B_2=\left[ {\begin{array}{ccc}
3 & -2 & 1  \\
2 & 1 & 3 \\ 
-1 & -3 & 2  \\
\end{array} } \right]
$$
$c_1=c_2=[0, 0, 0]^t$ and $K= [-10,10]\times[-10,10]\times[-10,10]$.
  We performed some numerical experiments for this example, and hence we chose seven starting points. Our stopping criterion is
$\|v^{k-1}-v^k\|<\varepsilon$, and we take $\varepsilon=10^{-6}$.
 
The numerical results are displayed in the following table, where  the starting points, the obtained solution, the number of iterations and the  CPU time have been reported.

Also, for each starting point, the test was successful, meaning that the sequence $\{v^k\}$  converges to a solution of VQEP$(f,T)$. All problems were solved by the Optimization Toolbox in
Matlab R2020a on a Laptop Intel(R) Core(TM) i7-
8665U CPU @ 1.90GHz RAM 8.00 GB.

\begin{center}

{\scriptsize
\begin{tabular}{|p{2.5cm}|p{3cm}| p{3cm}| p{3cm}|}
 \hline
 \multicolumn{4}{|c|}{Experiment for Example \ref{ex-matrix} } \\
 \hline
 Starting point: $v^0$ & solution found &  Number of iterations &  CPU time (Sec) \\
 \hline
  (4, 2, -3)& (10, 10, 10) & 23 & 9.4218\\
  (5, -2, -5)& (10, 10, 10) & 69 & 13.4531\\
 (6, 3, -2)& (10, 10, 10) & 21 & 6.2968\\
 (-4, -3,-1)& (10, 10, 10) & 67 & 16.6718\\
 (4, 4, 4) & (-10, 10, 10) & 73 & 15.6875\\
 (7, -4, -3) & (-10, 10, 10) & 651 & 241.5937\\
 (5, -5, 5) & (-10, 10, 10) & 237 & 54.2343\\
  \hline
\end{tabular}
}
\end{center}
Note that both $(10, 10, 10)$ and $(-10, 10, 10)$ are solutions of the problem and  Theorem \ref{main-strong} says that the sequence $\{v^k\}$ generated by Algorithm SEML,
converges to a solution of the problem.
\end{example}

We end this paper by performing some numerical experiments.

\begin{example}\label{ex-4.5}
We define the vector valued bifunction
 $f:\mathbb{R}^2\times \mathbb{R}^2\to\mathbb{R}^2$ by
 \begin{equation}\label{bifun-ab}
f(x,y)=\Big(a(x_1^2+x_2)(y_1^2+y_2^2-x_1^2-x_2^2),  bx_1^2(y_1-x_1)+cx_2(y_2-x_2)\Big)
 \end{equation}
for all $x=(x_1, x_2), y=(y_1, y_2)\in \mathbb{R}^2$ 
where $a,b,c \in \mathbb{R}_+$. Consider   $C=\Big\{z\in \mathbb{R}^2: z_i\geq0, i=1,2\Big\}$ and $K=[-10,10]\times [1,10]$. 
We also define $T:K\to {\cal P}(K)$ by
\begin{equation}\label{maptx-1}
T(x)=\Big \{ z\in K: z_1+z_2\geq \max\{x_1+x_2, 2\}\Big\}.
 \end{equation}

It is obvious that  $f$ satisfies B1--B3 and  $T$ is a multivalued mapping with nonempty, closed and convex values. Also, similar to Example \ref{ex-matrix}, it can be shown that  $T$ satisfies B4 with $g(\cdot)=\|\cdot\|^2$. Note that $S(f,T)\not=\emptyset$; indeed, it is easy to check that for all $(a,b,c)\in\mathbb{ R}^3_{++}$  the unique solution is $x^*=(1,1)$.
 In order to implement our algorithm (SEML) in Section \ref{s3}, we  take $\delta=\frac{1}{1000}$, $\theta=\frac{1}{2}$,  $\beta_k\equiv 1$, $\gamma_k\equiv 1$ and $e^k\equiv (1, 1)$. If $\{v^k\}$ is the sequence generated by Algorithm SEML, then  Theorem \ref{main-strong} ensures that $\{v^k\}$
converges to the solution of VQEP$(f,T)$. We performed some numerical experiments for this example. We chose randomly $100$ random triples $(a,b,c)\in [0,100]\times[0,100]\times[0,100]$ and five starting points. Our stopping criterion is
$\|v^{k-1}-v^k\|<\varepsilon$, and we take $\varepsilon=10^{-6}$.
 
The numerical results are displayed in the following table, where  the starting points, the obtained solution, the average number of iterations and the average CPU times have been reported.

Also, all tests for the 100 problems corresponding to each starting point were successful, meaning that the sequence $\{v^k\}$  converges to $(1,1)$, which is  the solution of VQEP$(f,T)$. All problems were solved by the Optimization Toolbox in
Matlab R2020a on a Laptop Intel(R) Core(TM) i7-
8665U CPU @ 1.90GHz RAM 8.00 GB.

\begin{center}

{\scriptsize
\begin{tabular}{|p{2.5cm}|p{2.5cm}|p{4.5cm}| p{3.5cm}|}
 \hline
 \multicolumn{4}{|c|}{Experiment for Example \ref{ex-4.5} } \\
 \hline
 Starting point: $v^0$ & solution found &  Average number of iterations & Average CPU time (Sec) \\
 \hline
  (-3, 2)& (1, 1) & 78.33 & 13.1979\\
  (-9, 7)& (1, 1) & 64.84 & 14.1020\\
 (0, 2)& (1, 1) & 65.13 & 7.3579\\
 (2, 8)& (1, 1) & 34.24 & 7.6868\\
 (-5, 5)& (1, 1) & 64.24 & 9.7689\\
  \hline
\end{tabular}
}
\end{center}
\end{example}

\end{document}